\documentclass[a4paper,10pt]{article}
\usepackage{latexsym,amssymb,amsmath,amscd,amscd,amsthm,amsxtra}
\usepackage[dvips]{graphicx}

\newtheorem{thm}{Theorem}[section]
\newtheorem{lem}[thm]{Lemma}

\newtheorem{pro}[thm]{Proposition}
\newtheorem{ex}[thm]{Example}
\newtheorem{rmk}[thm]{Remark}
\newtheorem{defi}[thm]{Definition}

\newcommand{\lon }{\,\rightarrow\,}
\newcommand{\be }{\begin{eqnarray*}}
\newcommand{\ee }{\end{eqnarray*}}

\newcommand{\alphabetaarrow}[1]{\rightrightarrows{#1};\;\alpha,\beta}
\newcommand{\defbe}{\triangleq}
\newcommand{\pf}{\noindent{\bf Proof.}\ }

\newcommand{\inverse}{^{-1}}

\newcommand{\Real}{\mathbb R}

\newcommand{\huaB}{\huaV}

\newcommand{\huaA}{\huaU}

\newcommand{\huaF}{\mathcal{F}}

\newcommand{\huaU}{\mathcal{U}}
\newcommand{\huaV}{\mathcal{V}}

\newcommand{\huaC}{{\mathcal{C}}}

\newcommand{\huaT}{\mathcal{T}}

\newcommand{\CWM}{C^{\infty}(M)}
\newcommand{\CWN}{C^{\infty}(N)}

\newcommand{\XM}{\mathcal{X}(M)}
\newcommand{\XN}{\mathcal{X}(N)}

\newcommand{\set}[1]{\left\{#1\right\}}

\newcommand{\LeftMove}[1]{\overleftarrow{#1}}

\newcommand{\frkg}{\mathfrak g}
\newcommand{\frkh}{\mathfrak h}

\newcommand{\frkD}{\mathfrak D}

\newcommand{\frkI}{\mathfrak I}
\newcommand{\frkJ}{\mathfrak J}

\newcommand{\cawu}{\mathcal{A}}
\newcommand{\cawub}{\mathcal{B}}
\newcommand{\cawuc}{\mathcal{C}}

\newcommand{\otimesa}{{\otimes_{\cawu}}}
\def\gpd{\,\lower1pt\hbox{$\longrightarrow$}\hskip-.24in\raise2pt
         \hbox{$\longrightarrow$}\,}

\def\qed{\hfill ~\vrule height6pt width6pt depth0pt}

\newcommand{\wedgea}{\wedge_\cawu}
\newcommand{\wedgeb}{\wedge_\cawub}

\newcommand{\EXinga}{\ME^*_\cawu}
\newcommand{\FXingb}{{\MF}_{\cawub}^*}

\newcommand{\Dera}{\mathrm{Der}(\cawu)}

\newcommand{\pisharp}{\pi^\sharp}
\newcommand{\LieG}{\frkg}
\newcommand{\LieH}{\frkh}

\newcommand{\ME}{\mathbf{E}}
\newcommand{\MF}{\mathbf{F}}
\newcommand{\MG}{\mathbf{G}}
\newcommand{\MD}{\mathbf{D}}

\newcommand{\RightRight}[2]
{\stackrel{(#1,#2)}{\rightrightarrows}}
\newcommand{\RightLeft}[2]
{\stackrel{(#1,#2)}{\rightleftarrows}}
\newcommand{\LeftRight}[2]
{\stackrel{(#1,#2)}{\leftrightarrows}}

\newcommand{\Graph}{\mathcal{G}}
\newcommand{\EquivClass}{\overline}
\newcommand{\basemap}{\phi}

\newcommand{\PSLA}{\mathfrak{P}}
\newcommand{\CPSLA}{\mbox{\it co}\mathfrak{P}}
\newcommand{\MRSM}{\mathcal{M}}
\newcommand{\CMRSM}{\mbox{\it co}\mathcal{M}}
\newcommand{\Id}{\mathrm{Id}}
\newcommand{\Hom}{\mathrm{Hom}}
\newcommand{\Ker}{\mathrm{Ker}}
\newcommand{\pr}{\mathrm{pr}}
\begin{document}
\title{{On (Co-)morphisms of Lie Pseudoalgebras and Groupoids
\thanks{1) Research partially supported by NSF of China (19925105) and
China Postdoctoral Science Foundation (20060400017). 2) This article
appeared in \emph{Journal of Algebra}, 316(2007): 1-31.}}}
\author{ Z. Chen and Z.-J. Liu\\
Department of Mathematics and LMAM, Peking University,
Beijing 100871, China\\
          {\sf email: chenzhuo@math.pku.edu.cn,\quad liuzj@pku.edu.cn} }

\date{}


\maketitle

\begin{abstract}
We give a unified description of  morphisms and comorphisms of Lie
pseudoalgebras,   showing that the both types of morphisms can be
regarded as subalgebras of a  Lie pseudoalgebra, called the
$\psi$-sum. We   also provide similar descriptions for morphisms
and comorphisms of Lie algebroids and groupoids.
\end{abstract}
\section{Introduction}

As an algebraic version of Lie algebroids, a Lie pseudoalgebra is a
generalized structure for that of a Lie algebra over some field and
the structure of the $\cawu$-module $\Dera$, all derivations on  a
commutative algebra $\cawu$. It consists of a triple $(\ME, \cawu,~
\theta )$ where $\ME$ is an $\cawu$-module and a Lie algebra,
$\theta : \ME \longrightarrow \Dera$ (called the anchor of $\ME$),
which is  an $\cawu$-module morphism and a Lie algebra morphism
satisfying some compatibility conditions.

The language of Lie pseudoalgebras arises from that of Lie
algebroids, which were first introduced by Pradines \cite{Pradines}
to provide a precise description of the infinitesimal form of Lie
groupoids. The reader who wishes to pursue the topic of Lie
algebroids and groupoids is referred to Mackenzie's new book
\cite{Mkz:GTGA} (see also \cite[Chp \textsl{8}, \textsl{12}]{AKA},
\cite[I, III]{first}) for background information.
As an abstract algebraic treatment of the category of Lie algebroids
performed in \cite{PJHK} and \cite{Mackenzie:1995}, Lie
pseudoalgebras are variously called Lie Rinehart algebras
\cite{Rinehart,Huebschmann,Huebschmann2,Huebschmann3,Bkouche}, Lie
d-rings \cite{Palais}, Palais pairs \cite{G-S2}, differential Lie
algebras \cite{KSM}, or modules with differential \cite{PP}. Some
other closely related variants are Elie Cartan spaces
\cite{Barros1,Barros2}, Lie modules \cite{Nelson}, Lie-Cartan pairs
\cite{Kastler}(this list is not complete.) It may be seen as an
algebraic form of the notion of Lie algebroid in which vector
bundles over manifolds are replaced by modules over rings, vector
fields by derivations of rings, and so on.

We especially study morphisms and comorphisms, which determine two
different categories of Lie pseudoalgebras. The general notion of
morphisms of Lie algebroids comes from Higgins and Mackenzie
\cite{PJHK}. Comorphisms of Lie pseudoalgebras seem to have first
been defined by Huebschmann \cite{Huebschmann}. There are also
corresponding concepts of morphisms for Lie groupoids
\cite{PJHK2}. Morphisms of Lie pseudoalgebras are well known (for
example, see \cite{Huebschmann}). The analogue concept, i.e.,
comorphisms of groupoids can also be found in \cite{PJHK2}. Our
treatment of morphisms and comorphisms of Lie algebroids follows
from Higgins and Mackenzie \cite{PJHK2} and \cite{Mkz:GTGA}(see
also \cite{Mackenzie:1995}).

One of the main aim of this paper is to show that  the two kinds of
morphisms of Lie pseudoalgebras can be unified via restriction
theory. The conclusion is that both are subalgebras in a Lie
pseudoalgebra called the $\psi$-sum of two Lie pseudoalgebras. We
also develop similar theories for Lie algebroids as well as for Lie
groupoids. The reader will see that the three versions of our main
theorems, Theorems \ref{Thm:morphismLPA}, \ref{Thm:coandmorphismLAB}
and \ref{Thm:coandmorphismLGB}, are actually expressing the same
idea from different points of view.

The paper is organized as follows. In Section
\ref{Sec:Quotientrestriction}, we recall the basic definition of Lie
pseudoalgebras and we introduce the idea of restricting Lie
pseudoalgebras via an ideal of the algebra. Applying this process,
we obtain the so-called $\psi$-sum of two Lie pseudoalgebras, where
$\psi$ is an algebraic morphism.

Section \ref{Sec:MorandCophismofLPA} describes the two kinds of
morphisms of Lie pseudoalgebras. The principal objective in this
section is the proof of our  main result (Theorem
\ref{Thm:morphismLPA}) in this paper. It provides a picture of the
relationship of the two different morphisms, whose graphs turn out
to be two subalgebras of the $\psi$-sum with respect to a given
algebra morphism $\psi$. Moreover,
 as an application of  theorem \ref{Thm:morphismLPA},
we  give a short proof of the fact that with either morphisms or
comorphisms, the Lie pseudoalgebras form a category, which is
proved  originally in \cite{PJHK2}.

To gain a better understanding of these abstract theories, Section
\ref{Sec:MorandCophismofLAD} expresses the preceding theory in the
language of Lie algebroids. We recall the definition of morphism and
comorphism of Lie algebroids   and we show how both morphisms and
comorphisms are embedded as subalgebroids into a common algebroid
which we called the $\basemap$-sum. The $\basemap$-sum can be
regarded as the algebroid version of the preceding $\psi$-sum.

Section \ref{Sec:GeometricModel} is an exposition of the theory of
the two kinds of morphisms concerning groupoids, analogous to that
of Lie algebroids.  In this paper, we present an equivalent
definition of comorphism between groupoids with generally different
bases. It turns out that the action of a groupoid $\Gamma$ on
$Z\stackrel{\varphi}{\lon }M$ is actually a comorphism of groupoids
from $\Gamma$ to the pair groupoid $Z\times Z$.

As global versions of that of Lie algebroids, we then study the
relationship of comorphisms and morphisms of groupoids, with
likewise conclusions. Finally, we present various examples.

\vskip 0.2cm \noindent\textbf{Acknowledgements}

We acknowledge   the kind hospitality of ICTP where part of this
work was completed. We thank Y. Kosmann-Schwarzbach and A.
Weinstein for useful comments. Finally, special thanks to the
referee for many helpful suggestions and  pointing out typos and
erroneous statements in a previous version of this paper.

\section{The $\psi$-sum of Lie Pseudoalgebras}
\label{Sec:Quotientrestriction} Throughout the paper, $\cawu$ stands
for an associative,
 commutative algebra over $K$, where $K$ is the number
field\footnote{Our conclusions also hold if $K$ is assumed to be a
field of characteristic zero.} $\Real$ or $\mathbb{C}$. We shall be
concerned exclusively with algebras $\cawu$ which are
\emph{unitary}, and hence $K$ is regarded as a subset of $\cawu$. We
always consider left $\cawu$-modules, and they are denoted by
capital letters $\ME$, $\MF$, etc. By  ``$\cawu$-maps'', or ``maps
of $\cawu$-modules", we mean ``$\cawu$-module morphisms''. The dual
module of $\ME$, namely $\Hom_\cawu(\ME,\cawu)$, will be denoted by
$\EXinga$.

A derivation of $\cawu$ is a $K$-linear map $\delta$: $\cawu\lon
\cawu$, satisfying the Leibnitz rule
$\delta(ab)=\delta(a)b+a\delta(b)$, $\forall a,b\in\cawu$. The
$\cawu$-module $\Dera$ of all derivations of $\cawu$ is closed under
the bracket, or commutator $[\delta,\lambda]\defbe
\delta\circ\lambda-\lambda\circ\delta$ and $\Dera$ is a $K$-Lie
algebra.

In this section, we recall the basic definition of Lie
pseudoalgebras and   we mainly introduce the idea of restricting Lie
pseudoalgebras via an ideal of the algebra $\cawu$. Applying this
process, we obtain the so-called $\psi$-sum of two Lie
pseudoalgebras, where $\psi$ is an algebraic morphism.

\noindent\textbf{$\bullet$ Lie pseudoalgebras.}\\
 In this paper, we adopt
the following definition.
\begin{defi}\label{Def:DiffLieAlg}
Let $\ME$ be an $\cawu$-module and a $K$-Lie algebra with the
bracket $[~\cdot~,~\cdot~]$ : $\ME\times\ME\lon  \ME $. If there is
a  Lie algebra morphism $\theta$: $\ME\lon  \Dera$ {\rm (called the
anchor of $\ME$)} such that
\begin{eqnarray}\label{Eqt:bracketaxfy}
\theta(a_1  X_1)&=&a_1\theta( X_1),\\
\label{Eqt:bracketxfy} [ X_1,a_2 X_2]&=&a_2[ X_1, X_2]+\theta(
X_1)(a_2) X_2, \ \ \forall  X_i\in \ME, a_i\in \cawu,
\end{eqnarray}
we call $(\ME ,[~\cdot~,~\cdot~],\theta)$ a Lie pseudoalgebra over
$\cawu$.
\end{defi}
Abusing the notation of the brackets, we prefer to write
$$\theta(X)a=[X,a]=-[a,X],\ \ X\in \ME, a\in\cawu.$$
In fact, one can treat this convention as the Lie bracket in the
semidirect sum $\ME \ltimes \cawu$, where $\cawu$ is considered
 as an abelian Lie algebra.

In what follows, we write $(\ME ,\cawu)$ to indicate a Lie
pseudoalgebra $\ME$ over a commutative and associative $K$-algebra
$\cawu$, with the bracket and anchor $[~\cdot~,~\cdot~]_\ME$. And
for simplicity, we will sometimes suppress the subscripts
``${~}_\ME$".

When $\ME$ is a finitely generated projective $\cawu$-module, there
is an exterior differential operator $d_\ME: \wedgea^k \EXinga \lon
\wedgea^{k+1}\EXinga$ ($k\geq 0$), which is of square zero, defined
in the ordinary fashion\footnote{For finitely generated projective
$\cawu$-modules $\ME$, one has $\wedgea^k \EXinga\cong
\Hom_{\cawu}(\wedgea^k \ME,\cawu)$, and it is also finitely
generated projective. And $\ME\cong
\Hom_{\cawu}(\EXinga,\cawu)=(\EXinga)^*_{\cawu}$.}(see also
\cite{KSM,Huebschmann2,Huebschmann3}):
\begin{eqnarray}\nonumber
\langle d_\ME a,X\rangle &=&\theta(X)a=[X,a],\quad a\in\cawu,X\in
\ME;\\\nonumber \langle d_\ME \xi,X\wedgea Y\rangle
&=&\theta(X)\langle \xi,Y\rangle -\theta(Y)\langle \xi,X\rangle
-\langle \xi,[X,Y]\rangle ,\quad \xi\in \EXinga,X,Y\in
\ME\\\nonumber d_\ME(\xi_1\wedgea\cdots\wedgea \xi_n) &=&
(d_\ME\xi_1) \wedgea \xi_2\wedgea\cdots\wedgea \xi_n -\xi_1\wedgea
(d_\ME\xi_2) \wedgea\xi_3\wedgea\cdots\wedgea \xi_n
\\\label{Eqn:dDerivative} &&+\cdots+ (-1)^{n-1}
\xi_1\wedgea\cdots\wedgea \xi_{n-1}\wedgea (d_\ME\xi_n),
\quad\xi_i\in \EXinga.
\end{eqnarray}
Conversely, any operator $d_\ME$ satisfying the derivation rule
(\ref{Eqn:dDerivative}), and $d_\ME^2=0$, determines a Lie
pseudoalgebra structure on $\ME$ (over $\cawu$) \cite{KSM}.

\noindent\textbf{$\bullet$ The restriction of Lie pseudoalgebras}.\\
Fix a Lie pseudoalgebra $(\ME ,\cawu)$ and let $\frkI\varsubsetneq
\cawu$ be an ideal of $\cawu$. Then, we define
$$
 \ME^{\frkI}\defbe\set{X\in \ME| [X,\frkI]\subset \frkI},
$$
which is clearly a submodule of $\ME$, and $( \ME^{\frkI},\cawu)$
is a Lie pseudoalgebra. Let
$$\frkI \ME\defbe \set{\sum_i  a_iX_i| a_i\in \frkI, X_i\in \ME}.$$
Then $\frkI \ME\subset \ME^\frkI$ and the following lemma gives a
new Lie pseudoalgebra.
\begin{lem} The quotient
$(\cawu/\frkI)$-module $ \ME^{\frkI}/(\frkI \ME)$ is a Lie
pseudoalgebra with structures coming from $(\ME ,\cawu)$.
\end{lem}
\pf By definition, $\cawu (\frkI \ME)\subset \frkI \ME$, $\frkI
 \ME^{\frkI}\subset \frkI \ME$ and hence $ \ME^{\frkI}/(\frkI \ME)$ is
indeed an $(\cawu/\frkI)$-module. We define the induced brackets by
the obvious rules
$$[\EquivClass{X_1},\EquivClass{X_2}]\defbe
{[X_1,X_2]}^-,\quad [\EquivClass{X},\EquivClass{a}]\defbe
{[X,a]}^-.$$ Here by $\EquivClass{X}$ we denote the element $X+\frkI
\ME\in$ $ \ME^{\frkI}/(\frkI \ME)$, and similarly by
$\EquivClass{a}$ we denote $a+\frkI\in \cawu/\frkI$. It suffices to
prove that they are well defined. In fact, we have
$$[\frkI \ME, \cawu]\subset \frkI, \quad [ \ME^{\frkI},\frkI]\subset
\frkI, \quad [ \ME^{\frkI}, \frkI \ME]\subset \frkI \ME.$$ The first
two relations are   obvious. For the last one, notice that if
$X_1\in \ME^\frkI$, $X'=aX_0\in \frkI \ME$, where $a\in \frkI$,
$X_0\in \ME$, we have
$$[X_1,X']=a[X_1,X_0]-[X_1,a]X_0\in \frkI \ME.
\ \quad \quad \blacksquare$$

Thus we see that $( \ME^{\frkI}/(\frkI \ME),\cawu/\frkI)$ inherits
all the structures of $(\ME ,\cawu)$. As an application, in Example
\ref{Ex:QuotientReduction}, one will find how a Lie algebroid is
restricted to a Lie pseudoalgebra through a submanifold. But
generally speaking, we can only depict the quotient module $
\ME^{\frkI}/(\frkI \ME)\subset \ME/(\frkI \ME)$ in the following
manner.
\begin{lem}\label{Lem:EoverIEcongEotimesaB}
~~~~

\begin{itemize} \item[1)] Let $\ME$ be an $\cawu$-module and
$\frkI$ be an ideal of $\cawu$. Then, $\ME/(\frkI \ME)\cong
\ME\otimesa (\cawu/\frkI)$ as $(\cawu/\frkI)$-modules, under the
isomorphism $\sigma:$ $\bar{X}\mapsto  X\otimesa \bar{1}$.

\item[2)] Let $(\ME ,\cawu)$ be a Lie pseudoalgebra and $\frkI$ be
an ideal of $\cawu$. Then under the isomorphism $\sigma$ defined
above, the quotient module $ \ME^{\frkI}/(\frkI \ME)\cong
\ME^\frkI\otimesa (\cawu/\frkI)$, and the latter  has the induced
Lie pseudoalgebra structure {\rm(over $\cawu/\frkI$)} defined by
\begin{eqnarray}\label{Eqn:brackets1}
{[X_1 \otimesa \EquivClass{a_1},\EquivClass{a_2}]}
&=&\EquivClass{a_1} {[X_1,a_2]}^-,\\\label{Eqn:brackets2} {[X_1
\otimesa \EquivClass{a_1}, X_2 \otimesa \EquivClass{a_2}]}
&=&[X_1,X_2] \otimesa \EquivClass{a_1a_2}+ X_2 \otimesa
\EquivClass{a_1} {[X_1,a_2]}^-- X_1 \otimesa \EquivClass{a_2}
{[X_2,a_1]}^-.
\end{eqnarray}
for all $X_i\in \ME^\frkI$, $a_i\in\cawu$.
\end{itemize}
\end{lem}
\pf For 1), consider the map $\chi: \ME \times \cawu/\frkI \lon
\ME/(\frkI \ME)$, $(X,\bar{a})\mapsto \EquivClass{aX}$, $\forall
X\in \ME$, $a\in \cawu$. It follows that
$$\chi(a_1X,\EquivClass{a_2})=\EquivClass{a_1a_2X}=
a_1\chi(X,\EquivClass{a_2})=\chi(X,a_1\EquivClass{a_2}),$$ i.e.,
$\chi$ is $\cawu$-bilinear and hence it induces a well-defined
$\cawu$-map $\tilde{\chi}: \ME\otimesa (\cawu/\frkI)\lon \ME/(\frkI
\ME)$, $X\otimesa\EquivClass{a}\mapsto \EquivClass{aX}$. One can
check that $\chi$ is exactly the inverse map of $\sigma$.

Now under $\sigma$, the image of $ \ME^{\frkI} /(\frkI \ME)$ is
clearly $ \ME^{\frkI}\otimesa (\cawu/\frkI)$, and the bracket in $
\ME^{\frkI}/(\frkI \ME)$ is then transferred to $
\ME^{\frkI}\otimesa (\cawu/\frkI)$, simply given by
$$[X_1\otimesa\bar{1},\EquivClass{a}]={[X_1,a]}^-,\quad
[X_1\otimesa\bar{1}, X_2\otimesa\bar{1}]=[X_1,X_2]\otimesa
\bar{1}.$$ Equalities (\ref{Eqn:brackets1}) and
(\ref{Eqn:brackets2})  then follow   the bracket laws
(\ref{Eqt:bracketaxfy}) and (\ref{Eqt:bracketxfy}). \qed

This fact says that  $ \ME^{\frkI}/(\frkI \ME)$ can be regarded as
an $(\cawu/\frkI)$-submodule of $\ME \otimesa (\cawu/\frkI)$. But
the latter one does NOT have a Lie pseudoalgebra structure.

\begin{defi}
We   denote $$ \ME _\frkI=\ME^{\frkI}\otimesa
(\cawu/\frkI)=\ME^{\frkI}/(\frkI \ME)$$  and call $(\ME
_{\frkI},\cawu/\frkI)$ the \textbf{$\frkI$-restriction} of a Lie
pseudoalgebra $(\ME ,\cawu)$ with respect to the ideal $\frkI\subset
\cawu$.
\end{defi}
\begin{rmk}\rm The terminology ``restriction'' here is suggested by the referee from the point of view of geometry as follows.
Consider the $\CWM$-module $\XM$, the set of smooth vector fields on
a smooth manifold $M$. Let $N\subset M$ be a submanifold and
consider an ideal $\frkI\subset \CWM$ which is the collection of the
functions that vanish  on $N$. In this case, it is clear that
$\CWM/\frkI\cong\CWN$ and   $(\XM)_{\frkI}=\XN$. So, the
$\frkI$-restriction of $\XM$ is the usual restriction of the tangent
vector fields on a submanifold. The construction of such
restrictions in the Lie algebroid case  also corresponds to
restricting a Lie algebroid to a submanifold of the base, as Example
\ref{Ex:QuotientReduction} later makes clear.
\end{rmk}
\noindent\textbf{$\bullet$ The $\psi$-sum}. \\  A useful
application of the preceding restriction theory is as follows.
Consider two Lie pseudoalgebras $(\ME ,\cawu$) and $(\MF ,\cawub)$
and define their \textbf{direct sum} to be the
$\cawu\otimes\cawub$-module $(\ME \otimes\cawub)\oplus
(\cawu\otimes \MF)$. Here ``$\otimes$'' means ``$\otimes_K$'' and
we adopt this convention throughout the paper. We first endow the
direct sum with a new Lie pseudoalgebra structure out of the given
ones.
\begin{pro} $((\ME \otimes\cawub)\oplus
(\cawu\otimes \MF), \cawu\otimes\cawub)$ is a Lie pseudoalgebra, the
structure being given by the following rules:
\begin{eqnarray*}
&&[X_1\otimes b_1+a_1\otimes Y_1,a_2\otimes b_2]\\
&=&[X_1,a_2]_\ME\otimes b_1b_2 +a_1a_2\otimes [Y_1,b_2]_\MF\,;\\
&&[X_1\otimes b_1+a_1\otimes Y_1,X_2\otimes b_2+a_2\otimes Y_2]\\
&=&[X_1,X_2]_\ME\otimes b_1b_2 + a_1a_2\otimes [Y_1,Y_2]_\MF \\
&&+ [X_1,a_2]_\ME\otimes b_1Y_2 +[a_1,X_2]_\ME\otimes b_2Y_1 +
a_1X_2\otimes [Y_1,b_2]_\MF + a_2X_1\otimes [b_1,Y_2]_\MF.
\end{eqnarray*}
Here $a_1,a_2\in\cawu$, $b_1,b_2\in\cawub$, $X_1,X_2\in \ME$,
$Y_1,Y_2\in \MF$.
\end{pro}

We omit the proof of this proposition, since the computations are
quite straightforward. For simplicity, in what follows, we will
again omit the subscripts ``${~}_\ME$, ${~}_\MF$'', etc. One is also
referred to Example \ref{Ex:DirectSum} of the direct sum of two Lie
algebroids as a guide to understanding $((\ME \otimes\cawub)\oplus
(\cawu\otimes \MF),\cawu\otimes\cawub)$.

From now on, we fix an algebra morphism $\psi: \cawu\lon \cawub$.
And hence $\cawub$-modules can be regarded as $\cawu$-modules via
$\psi$. By $\tilde{\psi}$ we denote the morphism given canonically
by
\begin{equation}\label{Eqt:Psi}
\cawu\otimes\cawub\lon \cawub:\quad a\otimes b\mapsto \psi(a)b.
\end{equation}
Notice that, $\tilde{\psi}$ is always surjective and  $\cawub\cong
(\cawu\otimes\cawub)/\Ker\tilde{\psi}$. It is easy to see that
\begin{equation}\label{Eqt:kerPsi}
\frkJ\defbe \Ker\tilde{\psi}=Span_K\set{a\otimes b-1\otimes
\psi(a)b| a\in\cawu,b\in\cawub}.
\end{equation}

\begin{defi}Let $(\ME ,\cawu)$ and $(\MF ,\cawub)$ be two Lie
pseudoalgebras and let  $\psi: \cawu\lon \cawub$ be an algebraic
morphism. For the direct sum $\MD =(\ME \otimes\cawub)\oplus
(\cawu\otimes \MF)$ and the ideal $\frkJ=\Ker\tilde{\psi}$ defined
as above, we denote $(\MD _\frkJ,\cawub)$, the $\frkJ$-restriction
of $(\MD ,\cawu\otimes\cawub)$ by $(\ME \oplus_\psi \MF, \cawub )$
and call this pseudoalgebra the \textbf{$\psi$-sum} of $(\ME
,\cawu)$ and $(\MF ,\cawub)$ with respect to the morphism $\psi$.
\end{defi}
The following theorem gives an explicit description of $\ME
\oplus_\psi \MF $.

\begin{thm}\label{Thm:EopluspsiF}With the above notations,
\begin{itemize}
\item[1)] the $\psi$-sum
$\ME \oplus_\psi \MF$ is a $\cawub$-submodule of $(
\ME\otimesa\cawub)\oplus \MF$ {\rm (here the tensor product $\ME
\otimesa\cawub$ is both a $\cawub$ and $\cawu$-module, considering
$\cawub$ as an $\cawu$-module through $\psi$)};
\item[2)]
an element
 $\sum_i X_i\otimesa b_i+Y\in (
\ME\otimesa\cawub)\oplus \MF$ belongs to $\ME\oplus_\psi \MF$ if and
only if
\begin{equation}\label{Eqt:Descriptionofpsisum}
\sum_i \psi([X_i,a])b_i=[Y,\psi(a)], \ \forall a\in \cawu .
\end{equation}
\end{itemize}
\end{thm}
\pf For 1), by definition, the $\psi$-sum
$$\ME \oplus_\psi \MF =\MD^\frkJ\otimes_{\cawu\otimes\cawub}\cawub\subset
\MD\otimes_{\cawu\otimes\cawub}\cawub\cong( \ME\otimesa\cawub)\oplus
\MF.$$

For 2), since $\MD^\frkJ=\set{x\in \MD | \tilde{\psi}[x,\frkJ]=0}$,
and according to (\ref{Eqt:kerPsi}), an element $\sum_i (X_i\otimes
b_i+a_i\otimes Y_i)\in \MD ^\frkJ$ if and only if
$$ \tilde{\psi}(\sum_i[X_i\otimes b_i+a_i\otimes Y_i,
a\otimes b-1\otimes \psi(a)b]) =
(\sum_i\psi([X_i,a])b_i-[\sum_i\psi(a_i)Y_i,\psi(a)])b=0,$$ holds
for all $a\in\cawu$, $b\in\cawub$. This expression shows that
$$
\MD^\frkJ=\set{X_i\otimes b_i+a_i\otimes
Y_i|\sum_i\psi([X_i,a])b_i-[\sum_i\psi(a_i)Y_i,\psi(a)]=0,\forall
a\in\cawu}.
$$
On the other hand,
$$\sum_i (X_i\otimes b_i+a_i\otimes Y_i)\otimes_{\cawu\otimes\cawub}\cawub
=\sum_i X_i\otimesa b_i + \sum_i \psi(a_i)Y_i\,.$$ Thus, $\ME
\oplus_\psi \MF =\MD^\frkJ\otimes_{\cawu\otimes\cawub}\cawub$ is
exactly the $\cawub$-submodule described in this theorem.\qed

By Lemma \ref{Lem:EoverIEcongEotimesaB}, the reader should have no
difficulty in obtaining the expressions of brackets of $\ME
\oplus_\psi \MF $, given in the following proposition.
\begin{pro}\label{Pro:StructureEpluspsiF}
The structure maps of the Lie pseudoalgebra $(\ME \oplus_\psi \MF
,\cawub)$ are given as follows.
\begin{eqnarray}\label{Eqt:AnchorEpluspsiF}
[\sum_i X_i\otimesa b_i+Y,b]&=&[Y,b];\\\nonumber [\sum_i X_i\otimesa
b_i+Y,\sum_j X'_j\otimesa b'_j+Y'] &=& \sum_{i,j} [X_i,X'_j]\otimesa
b_ib'_j+\sum_j X'_j\otimesa [Y,b'_j]\\\nonumber &&- \sum_i
X_i\otimesa [Y',b_i] +[Y,Y'],\\\label{Eqt:BracketEpluspsiF}
\end{eqnarray}
for all $\sum_i X_i\otimesa b_i+Y$, $\sum_j X'_j\otimesa
b'_j+Y'\in \ME\oplus_\psi \MF $, $b\in\cawub$. \qed
\end{pro}

\noindent\textbf{$\bullet$  Examples.}\\ We will denote a Lie
algebroid $\huaA$ over a manifold $M$ by $(\huaA, M, \rho_\huaA )$,
or simply $(\huaA,M)$, with the anchor map $\rho_\huaA: \huaA\lon
TM$ and the {\rm($\Real$-)} Lie algebra structure
$[~\cdot~,~\cdot~]_\huaA$ on $\Gamma(\huaA)$. We also write
$[v,~\cdot~]_\huaA$ to indicate the tangent vector $\rho_\huaA(v)$,
for $v\in\huaA$. The differential operator on $\Gamma(\wedge^\bullet
\huaA^*)$ will be denoted by $d_\huaA$. When it will introduce no
confusion, these subscripts ``$_\huaA$" can be omitted.

\begin{ex}\label{Ex:QuotientReduction}
\rm The reader should bear in mind this example as a guide to what
is going on for the restriction theories. Let $(\huaA,N,\rho)$ be a
Lie algebroid, and let ${N_0}$ be an embedded submanifold. Let
$\iota: {N_0}\lon N$ be the inclusion. Consider $\cawu=\CWM$. Of
course $\frkI=\Ker \iota^*=\set{f\in\cawu| f|_{N_0}=0}$ is an ideal
of $\cawu$. Denoting $\ME =\Gamma(\huaA)$,  we want to find the
restriction of the Lie pseudoalgebra $(\ME ,\cawu)$ with respect to
$\frkI$. Clearly,
$$ \ME^{\frkI}=\set{A\in \Gamma(\huaA)| \rho(A)|_{N_0}\in T{N_0}},$$
$$\frkI \ME=\set{A\in \Gamma(\huaA)| A|_{N_0}=0}.$$
Thus, the quotient algebra $\ME _\frkI$ can be regarded as the
 space of sections $$\huaA_{N_0}\defbe \set{v\in \huaA_p| p\in {N_0},
\rho(v)\in T_p{N_0}},$$ as an $\cawu/\frkI\cong
C^\infty{({N_0})}$-module.

Although $\huaA_{N_0}$ may not be a vector bundle over ${N_0}$
(since the dimension of its fibers may vary), it has formally the
anchor map and   the Lie algebra structure of the  space of sections
$\Gamma(\huaA_{N_0})$. Thus $\huaA_{N_0}$ can be regarded
as a \textit{generalized} Lie algebroid. 
\end{ex}
\begin{ex}\rm
Let $\ME$ and $\MF$ be two $K$-Lie algebras. The only $K$-morphism
from $K$ to $K$ is the identity map $I$, and the $I$-sum of
$(\ME,K)$ and $(\MF,K)$ is just $(\ME\oplus \MF,K)$, since in this
case (\ref{Eqt:Descriptionofpsisum}) becomes $0=0$.
\end{ex}
\begin{ex}\rm
Let $(\ME,\cawu)$ be a Lie pseudoalgebra and $\MF$ be a $K$-Lie
algebra. For an algebra morphism $\psi: \cawu\lon K$, the
$\psi$-sum of $(\ME,\cawu)$ and $(\MF,K)$ is the set
$$\ME\oplus_\psi \MF =
\set{ X\otimesa 1+ Y| X\in \ME,Y\in \MF, [X,\cawu]_\ME\subset
\Ker\psi}.$$ For example, when $\ME=\Gamma(\huaA)$, $\cawu=\CWN$
where $(\huaA,M)$ is a Lie algebroid, and $\psi(h)=h(p)$, $\forall
h\in \CWN$, for some fixed point $p\in N$, we have
$$\ME\oplus_\psi \MF
=\set{(v,Y)\in \huaA_p\times \MF| \rho_\huaA(v)=0}.$$
\end{ex}
\begin{ex}\label{Ex:DirectSum}\rm
We recall the direct product of Lie algebroids as in Higgins and
Mackenzie \cite{PJHK} and Mackenzie \cite{Mkz:GTGA}. Let $(\huaA, M,
\rho_\huaA )$ and $(\huaB, N, \rho_\huaB )$ be two Lie algebroids.
Their direct sum is a bundle $\huaA\times \huaB$ over $M\times N$,
with elements $(u,v)$, $u\in \huaA_p$, $v\in \huaB_q$, bundle map
$(u,v)\mapsto (p,q)$. We make $\huaA\times \huaB$ into a Lie
algebroid on the base $M\times N$ as follows. The anchor map is
$(u,v)\mapsto (\rho_\huaA(u),\rho_\huaB(v))$. For two sections of
$\huaA\times \huaB$, we define the Lie bracket to be
\begin{eqnarray*}
&&[(f_1A_1,g_1B_1),(f_2A_2,g_2B_2)]\\
&=&(f_1f_2[A_1,A_2]_\huaA+g_1[B_1,f_2]_\huaB A_2
-g_2[B_2,f_1]_\huaB A_1,\\&& \quad
g_1g_2[B_1,B_2]_\huaB+f_1[A_1,g_2]_\huaA B_2 -f_2[A_2,g_1]_\huaA
B_1).
\end{eqnarray*}
Here $f_i\in \CWN$, $g_i\in \CWM$, $A_i\in \Gamma(\huaA)$, $B_i\in
\Gamma(\huaB)$. To see that $(\Gamma(\huaA\times
\huaB),C^\infty(M\times N))$ is exactly the direct sum of
$(\Gamma(\huaA),\CWM)$ and $(\Gamma(\huaB),\CWN)$, one should regard
$C^\infty(M\times N)\thickapprox \CWM\otimes \CWN$ and
$\Gamma(\huaA\times \huaB)$ as the $C^\infty(M\times N)$-module
$\Gamma(\huaA)\otimes \CWN\oplus \Gamma(\huaB)\otimes \CWM$.


Given a smooth map $\basemap: M\lon N$, one has its graph
$$\Graph(\basemap)=\set{(x,\basemap(x))|x\in M}\subset M\times
N.$$ Applying the restriction process in Example
\ref{Ex:QuotientReduction} to the direct sum of Lie algebroids, we
define the $\basemap$-sum of Lie algebroids $\huaA$ and $\huaB$ to
be $(\huaA\times \huaB,M\times N)_{\Graph(\basemap)}$ (notice that
the result is a subset of $\huaA\oplus \basemap^!\huaB$), i.e.,

$$
\huaA\oplus_{\basemap}\huaB \defbe  \set{(u,v)\in \huaA\oplus
\basemap^!\huaB|\forall x\in M, u\in \huaA_x, v\in
\huaB_{\basemap(x)},\mbox{   }
{\basemap_*}\circ\rho_\huaA(u)=\rho_\huaB(v)}.
$$

The  space of sections $\Gamma(\huaA\oplus_{\basemap}\huaB)$ is a
Lie pseudoalgebra over $\CWM$, although
$\huaA\oplus_{\basemap}\huaB$ may not be a bundle over $M$. So it is
actually a generalized Lie algebroid,  an algebroid version
corresponding to the $\basemap^*$-sum of $(\Gamma(\huaA),\CWM)$ and
$(\Gamma(\huaB),\CWN)$.
\end{ex}

\section{Morphisms and Comorphisms of Lie Pseudoalgebras}\label{Sec:MorandCophismofLPA}
This section describes the two kinds of morphisms of Lie
pseudoalgebras. As a main result of this paper, theorem
\ref{Thm:morphismLPA} provides
 a picture of the relationship of the two different morphisms, whose
graphs turn out to be two subalgebras of the $\psi$-sum with respect
to a given algebra morphism $\psi$.

It is proved in \cite{PJHK2} that under either morphisms or
comorphisms, and together with certain composition laws, the Lie
pseudoalgebras are the objects of a category (in sense of morphism
or comorphism). However, it is not  easy to prove that the
comorphisms satisfy the usual axioms of categories. As an
application of  theorem \ref{Thm:morphismLPA}, we shall give a
short proof of the fact that with either morphisms or comorphisms,
the Lie pseudoalgebras form a category.

For an $\cawu$-module $\ME$ and a $\cawub$-module $\MF$, by  $(\ME
,\cawu)\RightRight{\Psi}{\psi}(\MF ,\cawub)$ we denote a pair of
morphisms $(\Psi, \psi)$, where $\psi: \cawu\lon \cawub$ is a
morphism of algebras and $\Psi: \ME\lon \MF$ is a map of
$\cawu$-modules {\rm{(}}considering $\cawub$-modules as
$\cawu$-modules through $\psi${\rm{)}}. We  call $\Psi$ an
$\cawu$-map over $\psi$.

Similarly, by writing $(\MF ,\cawub)\RightLeft{\Psi}{\psi}(\ME
,\cawu)$, we mean a pair of morphisms $(\Psi, \psi)$, where $\psi:
\cawu\lon \cawub$ is a morphism of algebras and $\Psi:\MF\lon
\ME\otimesa \cawub$ is a map of $\cawub$-modules. We also call
$\Psi$ a $\cawub$-map over $\psi$.

The concept of morphisms of Lie pseudoalgebras is given by
\cite{Huebschmann,Mackenzie:1995} as follows, and thus Lie
pseudoalgebras with their morphisms form a category.
\begin{defi}\label{Def:MorphismLPA}
Let $(\ME ,\cawu)$, $(\MF ,\cawub)$ be two Lie pseudoalgebras. A
\textbf{morphism} of Lie pseudoalgebras from $\ME$ to $\MF$, is a
pair of morphisms $(\ME ,\cawu)\RightRight{\Psi}{\psi}(\MF ,\cawub)$
such that
\begin{itemize}
\item[1)]
 $\psi([X,a])=[\Psi(X),\psi(a)],\quad\quad \forall X\in \ME,
a\in\cawu;$
\item[2)]
$\Psi([X_1,X_2])=[\Psi( X_1),\Psi( X_2)], \quad\quad\forall X_1,
X_2\in \ME.$
\end{itemize}
In particular, if both $\Psi$ and $\psi$ are injective, we call
$(\ME,\cawu)$ a \textbf{Lie subpseudoalgebra} of $(\MF,\cawub)$.
\end{defi}
Meanwhile, there is another kind of morphism defined for Lie
pseudoalgebras, also allowing the bases to be changed.

\begin{defi}\cite{Mackenzie:1995,PJHK2}\label{Def:CoMorphismLPA}
Let $(\ME ,\cawu)$, $(\MF ,\cawub)$ be two Lie pseudoalgebras. A
\textbf{comorphism} of Lie pseudoalgebras from $\MF$ to $\ME$ over
$\psi$, is a pair of morphisms $(\MF
,\cawub)\RightLeft{\Psi}{\psi}(\ME ,\cawu)$, such that
\begin{itemize}
\item[1)] if $Y\in \MF$ and $\Psi(Y)=\sum_k X_k\otimesa b_k$, for
some $X_k\in \ME$ and $b_k\in \cawub$, then
\begin{equation}\label{Eqt:CoMorphsim1}
[Y,\psi(a)]=\sum_k b_k\psi([X_k,a]),\quad \forall a\in\cawu;
\end{equation}

\item[2)] if $Y_1,Y_2\in \MF$ and $\Psi(Y_1)=\sum_k
X_{1,k}\otimesa b_{1,k}$, $\Psi(Y_2)=\sum_l X_{2,l}\otimesa
b_{2,l}$, where $X_{i,\cdot}\in \ME$, $b_{i,\cdot}\in \cawub$, then
\begin{equation}\label{Eqt:CoMorphsim2}
\Psi([Y_1,Y_2])=\sum_{k,l}[X_{1,k},X_{2,l}]\otimesa
b_{1,k}b_{2,l}+\sum_l X_{2,l}\otimesa [Y_1,b_{2,l}]-\sum_k
X_{1,k}\otimesa [Y_2,b_{1,k}]\,.
\end{equation}
\end{itemize}
In particular, if $\Psi$ is injective and $\psi$ is surjective, we
call $(\MF,\cawu)$ a \textbf{Lie co-subpseudoalgebra} of
$(\ME,\cawu)$.
\end{defi}
\begin{rmk}\rm~~~~
\begin{itemize}
\item[1)] Notice that $\psi$ here is in the reverse direction to that
in the definition of morphisms, and $\Psi$ is not a map from $\MF$
to $\ME$, but a map $ \MF$ to $\ME\otimesa \cawub$.

\item[2)] In the particular case that $\ME$ and $\MF$ are Lie
pseudoalgebras over the same algebra $\cawu=\cawub$, the  two
notions of morphisms are the same. Moreover, a Lie subpseudoalgebra
$(\ME,\cawu)$ of $(\MF,\cawu)$ is of course a Lie
co-subpseudoalgebra, and vice-versa.

\end{itemize}
\end{rmk}

The preceding Definition \ref{Def:MorphismLPA} seems very natural,
compared to that of the comorphism of Definition
\ref{Def:CoMorphismLPA}. The definition of comorphisms of Lie
pseudoalgebras appeared in \cite{Mackenzie:1995,PJHK2}. In \cite{PP}
Popescu had given similar concepts for various kinds of objects,
referred to modules with differentials. Although quite profound and
complicated, they arise from a clear geometric picture. The
morphisms of Lie pseudoalgebras correspond to the  morphisms of Lie
algebroids in Definition \ref{Def:MorphismLiealgebroidSimple} and
Theorem \ref{Thm:MorphismLiealgebroid}, which are the infinitesimal
form of the morphisms of Lie groupoids (see Definition
\ref{Def:MorphismGroupoid}).

If the modules are all finitely generated and projective, we have a
concise, yet equivalent description of comorphisms, from which the
reader will understand why the name ``comorphism" is used for such a
concept. First, we need the following proposition, concerning some
special properties of finitely generated projective modules.
\begin{pro}Let $\ME$, $\MF$ be finitely generated projective $\cawu$
and $\cawub$-modules respectively, and let $\psi: \cawu\lon
\cawub$ be an algebra morphism. Then
\begin{itemize}
\item[1)] $\ME \otimesa\cawub$ is a finitely generated projective
$\cawub$-module, and so are
$$(\ME \otimesa\cawub)\wedgea(\ME
\otimesa\cawub)= (\wedgea^2\ME) \otimesa\cawub, \quad\cdots ,
\quad(\wedgea^k \ME)\otimesa\cawub.$$

\item[2)] The map $I:$ $\ME \otimesa\cawub\lon
\Hom_{\cawu}(\EXinga,\cawub)$, sending each $X\otimesa b$ to
$$I( X\otimesa b) :\quad
\xi\mapsto \psi(\langle \xi,X\rangle )b,\quad\forall \xi\in \EXinga,
b\in\cawub,$$ is an isomorphism of $\cawub$-modules. Similarly, we
have
$$(\wedgea^2 \ME)\otimesa\cawub\cong
\Hom_{\cawu}(\wedgea^2\EXinga,\cawub),\,\cdots,\,(\wedgea^k
\ME)\otimesa\cawub\cong \Hom_{\cawu}(\wedgea^k\EXinga,\cawub).$$

\item[3)] Let $\Psi:\MF\lon  \ME\otimesa \cawub$ be a
$\cawub$-map. There is an induced $\cawu$-map $\Psi^*:\EXinga\lon
\FXingb$, called the dual map of $\Psi$, such that
$$\langle \Psi^*(\xi),Y\rangle =\langle I\circ\Psi(Y),\xi\rangle ,\quad\forall \xi\in \EXinga,Y\in \MF.$$

\item[4)] Let $\overline{\Psi}: \EXinga\lon \FXingb$ be a map of
$\cawu$-modules. There is a unique $\cawub$-map $\Psi:\MF\lon
\ME\otimesa \cawub$ such that $\Psi^*=\overline{\Psi}$, {\rm (
$\Psi(Y)=I\inverse\circ\langle\overline{\Psi}(\cdot),Y\rangle$, for
each $Y\in \MF$)}. \qed
\end{itemize}
\end{pro}
These four conclusions heavily rely on the condition that the
modules under consideration are finitely generated and projective.
We omit the proofs.
\begin{pro}\label{Pro:SimpleCoMorphismLPA}
Let $(\ME ,\cawu)$, $(\MF ,\cawub)$ be two finitely generated
projective Lie pseudoalgebras over the algebras $\cawu$ and
$\cawub$ respectively. Then the following two statements are
equivalent
\begin{itemize}
\item[1)] $(\MF ,\cawub)\RightLeft{\Psi}{\psi}(\ME ,\cawu)$ is a
comorphism of Lie pseudoalgebras;
\item[2)] the dual map $\Psi^*:\EXinga\lon \FXingb$ satisfies
\begin{equation}\label{Eqt:ChainMap}
d_\MF\circ \Psi^*=\Psi^*\circ d_\ME, \mbox{ as a map }
\wedgea^k\EXinga\lon \wedgeb^{k+1}\FXingb\ \ (k\geq 0).
\end{equation}

Here we regard $\Psi^*=\psi: \wedgea^0\EXinga=\cawu\lon
\wedgeb^0\FXingb=\cawub$, and $\Psi^*$ naturally lifts to an
$\cawu$-map $\wedgea^k\EXinga\lon \wedgea^k\FXingb$.
\end{itemize}
\end{pro}
\pf Due to formula (\ref{Eqn:dDerivative}),   the second statement
is equivalent to two conditions: for all $a\in\cawu$,
$d_\MF(\psi(a))=\Psi^*(d_\ME(a))$ and for all $\xi\in \EXinga$,
$d_\MF(\Psi^*(\xi))=\Psi^*(d_\ME(\xi))$. We prove that these two
conditions are equivalent to the first statement.

Suppose that $Y\in \MF$ and $\Psi(Y)=\sum_k X_k\otimesa b_k$, for
some $X_k\in \ME$ and $b_k\in \cawub$. Then by condition
$d_\MF(\psi(a))=\Psi^*(d_\ME(a))$, we get
\begin{eqnarray}\nonumber
[Y,\psi(a)]&=&\langle d_\MF(\psi(a)), Y \rangle =
\langle\Psi^*(d_\ME(a)),Y \rangle =\langle d_\ME(a),I\circ\Psi(Y)
\rangle\\\label{Eqn:temp1} &=& \sum_k \psi \langle X_k,d_\ME(a)
\rangle b_k = \sum_k b_k\psi([X_k,a]).
\end{eqnarray}
This proves relation (\ref{Eqt:CoMorphsim1}). Suppose that
$Y_1,Y_2\in \MF$ and $\Psi(Y_1)=\sum_k X_{1,k}\otimesa b_{1,k}$,
$\Psi(Y_2)=\sum_l X_{2,l}\otimesa b_{2,l}$, for some $X_{i,k}\in
\ME$, $b_{i,l}\in \cawub$. Then
\begin{eqnarray*}
&& \langle \Psi^*(d_\ME(\xi)), Y_1\wedgeb Y_2\rangle = \langle
d_\ME(\xi), I(\Psi(Y_1)\wedgeb \Psi(Y_2))\rangle\\
&=& \langle d_\ME(\xi), I(\sum_{k,l} (X_{1,k}\wedgea
X_{2,l})\otimesa b_{1,k}b_{2,l})\rangle\\
&=&\sum_{k,l} \psi\langle d_\ME(\xi), X_{1,k}\wedgea X_{2,l}
\rangle b_{1,k}b_{2,l}\\
&=&\sum_{k,l} \psi([X_{1,k},\langle \xi, X_{2,l} \rangle
]-[X_{2,l},\langle \xi, X_{1,k} \rangle ]-\langle \xi,
[X_{1,k},X_{2,l}]\rangle)b_{1,k}b_{2,l}\\
&=& \sum_l \langle d_\MF\circ \psi(\langle \xi, X_{2,l} \rangle) ,
Y_1\rangle b_{2,l}- \sum_k \langle d_\MF\circ \psi(\langle \xi,
X_{1,k} \rangle), Y_2 \rangle b_{1,k}- \sum_{k,l}\psi\langle \xi,
[X_{1,k},X_{2,l}]\rangle b_{1,k}b_{2,l}\,.
\end{eqnarray*}
Here we have used Equation (\ref{Eqn:temp1}). On the other hand,
we have
\begin{eqnarray*}
&& \langle d_\MF(\Psi^*\xi)), Y_1\wedgeb Y_2\rangle \\
&=& [Y_1,\langle\Psi^*\xi,Y_2\rangle]-
[Y_2,\langle\Psi^*\xi,Y_1\rangle]-\langle\Psi^*\xi,[Y_1,Y_2]\rangle\\
&=& [Y_1, \sum_l \psi\langle X_{2,l}, \xi \rangle b_{2,l}] -[Y_2,
\sum_k \psi\langle X_{1,k}, \xi \rangle b_{1,k}]
-\langle\xi,I\circ\Psi[Y_1,Y_2]\rangle
\\&=& \sum_l \langle
d_\MF(\psi\langle X_{2,l}, \xi \rangle b_{2,l}),  Y_1 \rangle
-\sum_k \langle d_\MF(\psi\langle X_{1,k}, \xi \rangle b_{1,k}),
Y_2\rangle-\langle\xi,I\circ\Psi[Y_1,Y_2]\rangle\\
&=& \sum_l \langle d_\MF\circ \psi(\langle \xi, X_{2,l} \rangle) ,
Y_1\rangle b_{2,l}- \sum_k \langle d_\MF\circ \psi(\langle \xi,
X_{1,k} \rangle), Y_2 \rangle b_{1,k}\\
&&+\sum_l \psi(\langle \xi, X_{2,l} \rangle) \langle d_\MF b_{2,l},
Y_1\rangle - \sum_k \psi(\langle \xi, X_{1,k} \rangle)\langle d_\MF
b_{1,k} , Y_2 \rangle-\langle\xi,I\circ\Psi[Y_1,Y_2]\rangle.
\end{eqnarray*} Therefore,
by condition $d_\MF(\Psi^*(\xi))=\Psi^*(d_\ME(\xi))$, we obtain
\begin{eqnarray}\nonumber
&&\langle\xi,I\circ\Psi[Y_1,Y_2]\rangle\\\label{Eqn:temp2} &=&
\sum_{k,l}\psi\langle \xi, [X_{1,k},X_{2,l}]\rangle
b_{1,k}b_{2,l}+\sum_l \psi(\langle \xi, X_{2,l} \rangle) \langle
d_\MF b_{2,l}, Y_1\rangle - \sum_k \psi(\langle \xi, X_{1,k}
\rangle)\langle d_\MF b_{1,k} , Y_2 \rangle.
\end{eqnarray}
Let $$Z=\sum_{k,l}[X_{1,k},X_{2,l}]\otimesa b_{1,k}b_{2,l}+\sum_l
X_{2,l}\otimesa [Y_1,b_{2,l}]-\sum_k X_{1,k}\otimesa
[Y_2,b_{1,k}].$$ Then, relation (\ref{Eqn:temp2}) says
$$\langle\xi, I(\Psi[Y_1,Y_2]-Z)\rangle=0.$$
By the arbitrariness of $\xi\in\EXinga$ and $I$ being an
isomorphism, we conclude that $\Psi[Y_1,Y_2]-Z=0$ and this proves
relation (\ref{Eqt:CoMorphsim2}). A similar process shows that 2)
implies 1).
 \qed

We  call $\Psi^*$ as a \textbf{chain map} if it enjoys the
property described by Equation (\ref{Eqt:ChainMap}).

\noindent\textbf{$\bullet$ The Main Theorem.}\\ We define the graph
of a pair of morphisms $(\ME ,\cawu)\RightRight{\Psi}{\psi}(\MF
,\cawub)$ where $\psi: \cawu\lon \cawub$ is a morphism of algebras
and $\Psi: \ME\lon \MF$ is a map of $\cawu$-modules.  Its graph is
defined  to be the $\cawub$-submodule
$$\Graph_{(\Psi,\psi)}\defbe \set{x+\widetilde{\Psi}(x)|x\in \ME\otimesa\cawub}
\subset( \ME\otimesa\cawub)\oplus\MF. $$  Here $\widetilde{\Psi}$ is
the $\cawub$-map $\ME \otimesa\cawub\lon \MF$ defined by the obvious
rule $X\otimesa b\mapsto \Psi(X)b$, $\forall X\in \ME$, $b\in
\cawub$.

For a pair of morphisms $(\MF ,\cawub)\RightLeft{\Psi}{\psi}(\ME
,\cawu)$, where $\psi: \cawu\lon \cawub$ is a morphism of algebras
and $\Psi:\MF\lon \ME\otimesa \cawub$ is a map of $\cawub$-modules,
we define its graph to be the $\cawub$-submodule
$$\Graph_{(\Psi,\psi)}\defbe \set{{\Psi}(Y)+Y|Y\in \MF}
\subset( \ME\otimesa\cawub)\oplus\MF. $$

Now, we state the main theorem. First, recall Theorem
\ref{Thm:EopluspsiF} which claims that the $\psi$-sum $\ME
\oplus_\psi \MF $ is a $\cawub$-submodule of $(
\ME\otimesa\cawub)\oplus\MF$. For the two pairs of morphisms above,
the following  theorem claims that they are a morphism or comorphism
of Lie pseudoalgebras if and only their graphs are contained in the
$\psi$-sum as Lie subpseudoalgebras.
\begin{thm}\label{Thm:morphismLPA}
Let $(\ME ,\cawu)$, $(\MF ,\cawub)$ be two Lie pseudoalgebras and
let $\psi: \cawu\lon \cawub$ be a morphism of algebras.
\begin{itemize}
\item[1)]Let
$(\ME ,\cawu)\RightRight{\Psi}{\psi}(\MF ,\cawub)$ be a pair of
morphisms. It is a morphism of Lie pseudoalgebras if and only if its
graph $\Graph_{(\Psi,\psi)}$ is a Lie subpseudoalgebra {\rm (over
$\cawub$)} of the $\psi$-sum $\ME \oplus_\psi \MF $.
\item[2)]Let
$(\MF ,\cawub)\RightLeft{\Psi}{\psi}(\ME ,\cawu)$ be a pair of
morphisms. It is a comorphism of Lie pseudoalgebras if and only if
its graph $\Graph_{(\Psi,\psi)}$ is a Lie subpseudoalgebra {\rm
(over $\cawub$)} of the $\psi$-sum $\ME \oplus_\psi \MF $.
\end{itemize}
\end{thm}
Before we prove these two statements, we state two facts which
deserve attentions.
\begin{lem}
~~~~~~
\begin{itemize}
\item[1)] With the same assumptions as in (1) of Theorem
\ref{Thm:morphismLPA}, $\Graph_{(\Psi,\psi)}$ is contained in $\ME
\oplus_\psi \MF $ if and only if (1) in Definition
\ref{Def:MorphismLPA} holds.

\item[2)] With the same assumptions as in (2) of Theorem
\ref{Thm:morphismLPA}, $\Graph_{(\Psi,\psi)}$ is contained in $\ME
\oplus_\psi \MF $ if and only if (1) in Definition
\ref{Def:CoMorphismLPA} holds.
\end{itemize}
\end{lem}
\pf 1) Consider $x= X\otimesa b$. Then $x+\widetilde{\Psi}(x)=
X\otimesa b+b\Psi(X)\in \Graph_{(\Psi,\psi)}$. By
(\ref{Eqt:Descriptionofpsisum}) in Theorem \ref{Thm:EopluspsiF},
$x+\widetilde{\Psi}(x)$ belongs to $\ME \oplus_\psi \MF $ if and
only if
$$\psi([X,a])b=[b\Psi(X),\psi(a)]=b[\Psi(X),\psi(a)]$$
holds. The conclusion then comes from the arbitrariness of $b$.

2) Let $Y\in \MF$ and suppose that $\Psi(Y)=\sum_k X_k\otimesa b_k$,
for some $X_k\in \ME$ and $b_k\in \cawub$. Using the equality
(\ref{Eqt:Descriptionofpsisum}) again, we know that
$\Psi(Y)+Y\in\Graph_{(\Psi,\psi)}$ belongs to $\ME \oplus_\psi \MF $
if and only if relation (\ref{Eqt:CoMorphsim1}) holds. \qed

\noindent\textbf{Proof of Theorem \ref{Thm:morphismLPA}}.

 1)
$\widetilde{\Psi}$ is clearly well defined. The above lemma says
that  (1) of Definition \ref{Def:MorphismLPA} is exactly the
condition that $\Graph_{(\Psi,\psi)}\subset \ME\oplus_\psi \MF $.
And under this condition, it is easily seen that (2) of Definition
\ref{Def:MorphismLPA} is exactly the condition for
$\Graph_{(\Psi,\psi)}$ to be closed under the bracket given by
relation (\ref{Eqt:BracketEpluspsiF}) in Proposition
\ref{Pro:StructureEpluspsiF}.


2) Let $Y_i\in \MF$ and suppose that $\Psi(Y_i)=\sum_j
X_{i,j}\otimesa b_{i,j}$, $i=1,2$. We have already shown that the
condition $\Psi(Y_i)+Y_i\in \ME\oplus_\psi \MF $ is  relation
(\ref{Eqt:CoMorphsim1}). And by Proposition
\ref{Pro:StructureEpluspsiF}, the condition that
$[\Psi(Y_1)+Y_1,\Psi(Y_2)+Y_2]\in \ME\oplus_\psi \MF $ is exactly
  relation (\ref{Eqt:CoMorphsim2}). \qed

It may  happen that, given $\psi$, no morphism or comorphism
$\Psi$ over $\psi$   exist. But the $\psi$- sum always exists.
This phenomenon also happens when Lie algebroids or groupoids are
concerned (see Remark \ref{Rmk:ItMayHappen} and  examples  at the
end of the paper).

\noindent\textbf{$\bullet$ The Category of Lie Pseudoalgebras.}\\
As an application of Theorem \ref{Thm:morphismLPA}, we now prove
that
 Lie pseudoalgebras are the objects of a category, with either
morphisms or comorphisms, which is proved  originally in
\cite{PJHK2}.  We will need the following fact which can be drawn
directly by the definition.
\begin{pro}\label{Pro:tripleSum}
Let $(\ME,\cawu)$, $(\MF,\cawub)$, $(\MG,\cawuc)$ be three Lie
pseudoalgebras and $\psi:\cawu\lon\cawub$, $\theta:\cawub\lon\cawuc$
be two algebra morphism. Then as a Lie subpseudoalgebra over
$\cawuc$,
$$(\ME\oplus_\psi \MF)\oplus_\theta \MG \subset \ME
\oplus_{\theta\circ\psi} (\MF\oplus_\theta \MG).$$
\end{pro}

The concept of category can be found in \cite{Hilton}. In the
present paper, we define the category $\PSLA$ and $\CPSLA$ of Lie
pseudoalgebras, using the following three pieces of data:
\begin{itemize}
\item[1)] a class of objects $(\ME,\cawu)$, $(\MF,\cawub)$,..., which are
respectively pseudoalgebras;
\item[2)] to each pair of objects $(\ME,\cawu)$, $(\MF,\cawub)$ of
$\PSLA$ (or $\CPSLA$), a set $\MRSM((\ME,\cawu);(\MF,\cawub))$, the
collection of all   morphisms $(\ME
,\cawu)\RightRight{\Psi}{\psi}(\MF ,\cawub)$ from $(\ME,\cawu)$ to
$(\MF,\cawub)$ (or, a set $\CMRSM((\ME,\cawu);(\MF,\cawub))$, the
collections of all   comorphisms $(\ME
,\cawu)\LeftRight{\Psi}{\psi}(\MF ,\cawub)$ from $(\MF,\cawub)$ to
$(\ME,\cawu)$);
\item[3)] to each triple of objects $(\ME,\cawu)$, $(\MF,\cawub)$,
$(\MG,\cawuc)$, a law of composition
\begin{equation}\label{Eqt:compositionMorphism}
\MRSM((\ME,\cawu);(\MF,\cawub))\times
\MRSM((\MF,\cawub);(\MG,\cawuc))\lon
\MRSM((\ME,\cawu);(\MG,\cawuc)).
\end{equation}
\begin{equation}\label{Eqt:compositionCoMorphism}
\mbox{(or}\quad \CMRSM((\ME,\cawu);(\MF,\cawub))\times
\CMRSM((\MF,\cawub);(\MG,\cawuc))\lon
\CMRSM((\ME,\cawu);(\MG,\cawuc)).\quad\mbox{)}\end{equation}
\end{itemize}
Here in (\ref{Eqt:compositionMorphism}), the composition of $(\ME
,\cawu)\RightRight{\Psi}{\psi}(\MF ,\cawub)$ and $(\MF
,\cawub)\RightRight{\Theta}{\theta}(\MG ,\cawuc)$ is the usual
composition of maps
$$(\ME ,\cawu)\RightRight{\Theta\circ\Psi}{\theta\circ\psi}(\MG
,\cawuc).$$ In (\ref{Eqt:compositionCoMorphism}), we define the
composition of $(\ME ,\cawu)\LeftRight{\Psi}{\psi}(\MF ,\cawub)$ and
$(\MF ,\cawub)\LeftRight{\Theta}{\theta}(\MG ,\cawuc)$ to be
$$(\ME ,\cawu)\LeftRight{\Psi*\Theta}{\theta\circ\psi}(\MG
,\cawuc).$$ Here $\Psi*\Theta: \MG\lon \ME\otimes_{\cawu}\cawuc$
(over $\theta\circ\psi$) is in fact $\overline{\Psi}\circ\Theta$,
where
$$\overline{\Psi}=\Psi\otimes \Id: \MF\otimes_{\cawub}\cawuc\lon
(\ME\otimesa\cawub)\otimes_{\cawub}\cawuc=\ME\otimesa\cawuc.$$ They
obviously satisfy the usual axioms of associativity and identity. So
the only problem is to check that the composition laws are well
defined. In what follows we prove that in
(\ref{Eqt:compositionCoMorphism}), the composition of two
comorphisms is still a comorphism. And one can do  similar procedure
for morphisms.

As above, let $(\ME ,\cawu)\LeftRight{\Psi}{\psi}(\MF ,\cawub)$ and
$(\MF ,\cawub)\LeftRight{\Theta}{\theta}(\MG ,\cawuc)$ be two
comorphisms. We prove that their composition $(\ME
,\cawu)\LeftRight{\Psi*\Theta}{\theta\circ\psi}(\MG ,\cawuc)$ is
also a comorphism of Lie pseudoalgebras. It suffices to show that
the graph
$$\Graph_{(\Psi*\Theta,\theta\circ\psi)}=
\set{{\Psi*\Theta}(Z)+Z|Z\in \MG} \subset(
\ME\otimesa\cawuc)\oplus\MG.$$ is a subalgebra of
$\ME\oplus_{\theta\circ\psi}\MG$ (over $\cawuc$). If we consider the
embedding map
$$
i: (\ME\otimes_{\cawu}\cawuc)\oplus \MG\lon
(\ME\otimes_{\cawu}\cawuc)\oplus (\MF\otimes_{\cawub}\cawuc)\oplus
\MG,
$$
$$
X\otimes_\cawu c + Z \mapsto X\otimes_\cawu c+ \Theta(Z)+Z,
$$
for all $X\in\ME$, $Z\in \MG$, $c\in\cawuc$, then it suffices to
show that $i(\Graph_{(\Psi*\Theta,\theta\circ \psi)})$ is a
subalgebra of
$$
i(\ME\oplus_{\theta\circ\psi}\MG)
=\ME\oplus_{\theta\circ\psi}\Graph_{(\Theta,\theta)}.
$$

In fact, by applying Theorem \ref{Thm:morphismLPA}, we know that
$\ME\oplus_{\theta\circ\psi}\Graph_{(\Theta,\theta)}$ is a
subalgebra of $\ME \oplus_{\theta\circ\psi} (\MF\oplus_\theta \MG)$.
And at the same time, $\Graph_{(\Psi,\psi)}\oplus_{\theta}\MG$ is a
subalgebra of $(\ME\oplus_{\psi}\MF)\oplus_{\theta}\MG$. So it is
also a subalgebra of $\ME \oplus_{\theta\circ\psi} (\MF\oplus_\theta
\MG)$ (by Proposition \ref{Pro:tripleSum}) and therefore the
intersection
$$
(\ME\oplus_{\theta\circ\psi}\Graph_{(\Theta,\theta)}) \cap
(\Graph_{(\Psi,\psi)}\oplus_{\theta}\MG)=i(\Graph_{(\Psi*\Theta,\theta\circ
\psi)})
$$
is a subalgebra of $\ME \oplus_{\theta\circ\psi} (\MF\oplus_\theta
\MG)$. Of course this implies that it is also a subalgebra of
$i(\ME\oplus_{\theta\circ\psi}\MG)$ and thus the proof is complete.

There are three crucial points in the above analysis:
\begin{itemize}
\item[1)] If $(\ME_1,\cawu)$ is a subalgebra of $(\ME_2,\cawu)$, and
$(\ME_2,\cawu)$ is a subalgebra of $(\ME_3,\cawu)$, then
$(\ME_1,\cawu)$ is a subalgebra of $(\ME_3,\cawu)$.
\item[2)] If $(\ME_1,\cawu)$ and $(\ME_2,\cawu)$ are both subalgebras
of $(\ME_3,\cawu)$, then their intersection $(\ME_1\cap
\ME_2,\cawu)$ is a subalgebra of $(\ME_3,\cawu)$.
\item[3)] If $(\ME_1,\cawu)$ is a subalgebra of $(\ME_3,\cawu)$, and
$(\ME_2,\cawu)$ is a subalgebra of $(\ME_3,\cawu)$, and
$\ME_1\subset \ME_2$, then $(\ME_1,\cawu)$ is a subalgebra of
$(\ME_2,\cawu)$.
\end{itemize}


\section{Morphisms and Comorphisms of Lie Algebroids}
\label{Sec:MorandCophismofLAD} This part is devoted to expressing
the preceding theories in the language of Lie algebroids. We recall
the definition of morphism and comorphism of Lie algebroids, which
originally appeared in \cite{PJHK2} and we recommend Mackenzie's
book \cite{Mkz:GTGA} for detailed information. However, we adopt a
different approach to these two concepts in this paper. Although the
original ones are equivalent to the definitions which follow, the
latter are quite concise and succinct in language. We finally show
how they are embedded into an algebroid which we called the
$\basemap$-sum, as subalgebroids.

A basic fact should be mentioned at the beginning. Consider a
smooth map $\basemap: M\lon N$. Let $V$ be a vector bundle over
$N$. We have the pull back bundle $\basemap^! V$ (over $M$) and a
morphism of algebras $\psi=\basemap^*: \CWN\lon \CWM$. For anther
vector bundle $W$ over $M$, and a bundle map $\Psi: \basemap^! V
\lon W$, it determines the dual bundle map
\begin{equation}
\begin{CD}
W^* @> \Psi^*>>   V^*\\
@V p VV @VV q V\\
M @> \basemap>>   N.
\end{CD}
\end{equation}
It also naturally induces an additive map $\widetilde{\Psi}$:
$\Gamma(V)\lon \Gamma(W)$ satisfying
\begin{equation}\label{Eqt:phipsiCompatible}
\widetilde{\Psi}(fB)=(\basemap^*f)\widetilde{\Psi}(B),\quad
\forall f\in\CWN, B\in\Gamma(V).
\end{equation}
This implies that $\widetilde{\Psi}$ is a map of $\CWN$-modules, in
the sense of $\basemap^*$. Conversely, any $\CWN$-module map
$\widetilde{\Psi}$: $\Gamma(V)\lon \Gamma(W)$ is uniquely determined
by such a bundle map $\Psi: \basemap^! V \lon W$, or $\Psi^*:
W^*\lon V^*$.

By definition, a Lie algebroid $(\huaA, M, \rho_\huaA )$ gives
rise to a Lie pseudoalgebra $(\Gamma(\huaA),\CWM)$ and we can
therefore extend the two kinds of morphisms described in Section
\ref{Sec:MorandCophismofLPA} to Lie algebroids, in both cases
allowing the bases to be
changed \cite{Mackenzie:1995,PJHK2,Mkz:GTGA}.

\begin{defi}\label{Def:SimpleCoMorphismLiealgebroid}
Let $(\huaA, M, \rho_\huaA )$ and $(\huaB, N, \rho_\huaB  )$ be
Lie algebroids on bases $M$ and $N$ respectively. Given a smooth
map $\basemap: M\lon N$, a \textbf{comorphism of Lie algebroids}
from $\huaB$ to $\huaA$ over $\basemap$ is a bundle map $\Psi:
\basemap^! \huaB \lon \huaA$, written
$(\huaB,N)\RightLeft{\Psi}{\basemap}(\huaA,M)$, such that the dual
map $\Psi^*: \huaA^*\lon\huaB^*$ is a Poisson map. Here $\huaA^*$
and $\huaB^*$ both carry the Lie-Poisson structures coming from
their Lie algebroid structures.

In particular, if $\basemap$ is surjective and $\Psi$ is  injective,
then we call $(\huaB,N)$ a \textbf{co-subalgebroid} of $(\huaA,M)$.
\end{defi}
According to the relationship between $\Psi$ and $\widetilde{\Psi}$,
we sometimes directly call $\widetilde{\Psi}: \Gamma(\huaB)\lon
\Gamma(\huaA)$ the comorphism of Lie algebroids. It will be
convenient at times to consider the comorphisms of Lie algebroids
from this alternate point of view, by using the following equivalent
description \cite{PJHK2}.
\begin{thm}\label{Thm:CoMorphismLiealgebroid}
With the assumptions in the above definition,
$(\huaB,N)\RightLeft{\Psi}{\basemap}(\huaA,M)$ is a comorphism of
Lie algebroids if and only if the following  two conditions hold:
\begin{itemize}\item[1)]
${\basemap_*}\circ\rho_\huaA\circ\Psi= \rho_\huaB$;

\item[2)] the induced map $\widetilde{\Psi}: \Gamma(\huaB)\lon
\Gamma(\huaA)$ is a morphism of Lie algebras.
\end{itemize}
\end{thm}
We remark that, the first condition can be restated as: for each
$B\in \Gamma(\huaB)$, the vector field
$\rho_\huaA(\widetilde{\Psi}(B))$ is $\basemap$-related to
$\rho_\huaB(B)$. The composition law of morphisms is
straightforward.

\begin{rmk}\rm Let $\psi=\basemap^*: \CWN\lon \CWM$.
In Theorem \ref{Thm:CoMorphismLiealgebroid}, (1) is equivalent to
$$\psi([B,g]_\huaB)=[\widetilde{\Psi}(B),\psi(g)]_\huaA,
\quad \forall B\in \Gamma(\huaB), g\in \CWN.$$ (C.f. relation (1) of
Definition \ref{Def:MorphismLPA}.) Hence the \emph{comorphism} in
Theorem \ref{Thm:CoMorphismLiealgebroid} is actually saying that
$$
(\Gamma(\huaB),\CWN)
\RightRight{\widetilde{\Psi}}{\psi}(\Gamma(\huaA),\CWM),
$$
is a \emph{morphism} of Lie pseudoalgebras.
\end{rmk}
\begin{ex}\rm
Let $(M,\pi)$  be a Poisson manifold and by $\Omega(M)$ we denote
$\Gamma(T^*M)$. With the $\pi$-bracket defined below
\begin{eqnarray*}
[\xi,f]_\pi &\defbe& \pisharp(\xi)(f);\\
\langle [\xi,\eta]_\pi,X\rangle  &\defbe &\langle
[X,\pi],\xi\wedge\eta\rangle +{\pisharp(\xi)}\langle \eta,X\rangle -
{\pisharp(\eta)}\langle \xi,X\rangle ,
\end{eqnarray*}
for any $\xi,\eta\in \Omega(M)$, $f\in\CWM$ and $X\in \XM$, it is
well known that $(\Omega(M),\CWM)$ is a Lie pseudoalgebra and $T^*M$
is a Lie algebroid on $M$. Let $(M,\pi)$ and $(N,\varpi)$ be two
Poisson manifolds. Assume that $\basemap: M\lon N$ is a Poisson map,
which induces ${\basemap_*}: TM\lon TN$ and the dual map
$\Phi=\basemap^*: \Omega(N)\lon \Omega(M)$. Then, $\Phi$ is a
comorphism of Lie algebroids over $\basemap$
$$(T^*N,N)\RightLeft{\basemap^*}{\basemap}(T^*M,M).$$
\end{ex}

The definition of morphisms of Lie algebroids appears in many texts
such as \cite{PJHK}, \cite{Mackenzie:1995} (see also
\cite{PJHK2,Mkz:GTGA}), stated in the form given in Theorem
\ref{Thm:MorphismLiealgebroid}. But we prefer to adopt a more
concise one as below.

\begin{defi}\label{Def:MorphismLiealgebroidSimple}
Let $(\huaA, M, \rho_\huaA )$ and $(\huaB, N, \rho_\huaB  )$ be
Lie algebroids on bases $M$ and $N$ respectively. A
\textbf{morphism of Lie algebroids} from $\huaA$ to $\huaB$,
written $(\huaA,M)\RightRight{\Psi}{\basemap}(\huaB,N)$, is a
vector bundle morphism
\begin{equation}
\begin{CD}
\huaA @> \Psi>>   \huaB\\
@V p VV @VV q V\\
M @> \basemap>>   N
\end{CD}
\end{equation}
such that the induced map $\widetilde{\Psi^*}:
\Gamma(\wedge^k\huaB^*)\lon \Gamma(\wedge^k\huaA^*)$ is a chain
map, i.e.,
$$ d_\huaA\circ
\widetilde{\Psi^*}=\widetilde{\Psi^*}\circ d_\huaB,\ \ \mbox{as a
map } \Gamma(\wedge^k\huaB^*)\lon \Gamma(\wedge^{k+1}\huaA^*) \ \
(k\geq 0).$$ Here we regard $\widetilde{\Psi^*}=\psi^*:
\Gamma(\wedge^0\huaB^*)=\CWN\lon \Gamma(\wedge^0\huaA^*)=\CWM$. In
particular, if $\basemap$   and $\Psi$ are both injective, then we
call $(\huaA,M)$ a \textbf{subalgebroid} of $(\huaB,N)$.
\end{defi}

\begin{thm}\label{Thm:MorphismLiealgebroid}
With the assumptions in the above definition,
$(\huaA,M)\RightRight{\Psi}{\basemap}(\huaB,N)$ is a morphism of Lie
algebroids if and only if
\begin{itemize}
\item[1)] \begin{equation}\label{Eqt:CoMorphsimLiealgebroids1}
\rho_\huaB\circ\Psi={\basemap_*}\circ \rho_\huaA; \end{equation}
\item[2)] if $A$, $A'\in\Gamma(\huaA)$ and ${\Psi^!}(A)=\sum_i f_i
B_i$, ${\Psi^!}(A')=\sum_i f'_j B'_j$, for some $f_i,f'_j\in
\CWM$, $B_i,B'_j\in\Gamma(\huaB)$, where ${\Psi^!}$ is the induced
bundle map $\huaA\lon \basemap^!\huaB$, then
\begin{equation}\label{Eqt:CoMorphsimLiealgebroids2}
{\Psi^!}([A,A'])=\sum_{i,j} f_if'_j[B_i,B'_j]+\sum_j [A,f'_j]_\huaA
B'_j-\sum_i [A',f_i]_\huaA B_i\,.
\end{equation}
\end{itemize}
\end{thm}
The reader should bear in mind that the second statement of the
theorem is a \emph{local} condition. The proof is omitted since it
is merely a repetition of Proposition \ref{Pro:SimpleCoMorphismLPA}.
Although the above two equalities seem quite complicated, by the
following remark and according to Proposition
\ref{Pro:SimpleCoMorphismLPA}, one is able to understand why
(\ref{Eqt:CoMorphsimLiealgebroids1}) together with
(\ref{Eqt:CoMorphsimLiealgebroids2}) are equivalent to Definition
\ref{Def:MorphismLiealgebroidSimple}.
\begin{rmk}\rm Let $\psi=\basemap^*: \CWN\lon \CWM$.
In Theorem \ref{Thm:MorphismLiealgebroid}, relation
(\ref{Eqt:CoMorphsimLiealgebroids1}) is equivalent to
$$[A,\psi(g)]_\huaA=
\sum_i f_i \psi([B_i,g]_\huaB),\quad \forall A\in \huaA, g\in
\CWN.$$ (C.f. relation (\ref{Eqt:CoMorphsim1}).) In addition, one
may regard
$\Gamma(\basemap^!\huaB)\thickapprox\Gamma(\huaB)\otimes_{\CWN}\CWM$.
Therefore, the \emph{morphism} in Theorem
\ref{Thm:MorphismLiealgebroid} is in fact a \emph{comorphism} of Lie
pseudoalgebras
$$
(\Gamma(\huaA),\CWM)\RightLeft{{\Psi^!}}{\psi}(\Gamma(\huaB),\CWN).
$$
\end{rmk}
Note that, by the two equivalent descriptions of morphisms of Lie
algebroids in Definition \ref{Def:MorphismLiealgebroidSimple} and
in Theorem \ref{Thm:MorphismLiealgebroid}, we recover Theorem 3.1
in \cite{HLZ}.

The above two kinds of morphisms are just the algebroid version
corresponding to those of Lie pseudoalgebras. Recall the
$\basemap$-sum $\huaA\oplus_{\basemap}\huaB$ of Lie algebroids
$\huaA$ and $\huaB$ define in Example \ref{Ex:DirectSum}. When
passing from the morphisms and comorphisms of Lie pseudoalgebras to
those of Lie algebroids, the relationships of the two kinds of
morphisms stated in Theorem \ref{Thm:morphismLPA}, admit
straightforward generalizations, formulated as follows.
\begin{thm}\label{Thm:coandmorphismLAB}
 Let $(\huaA, M, \rho_\huaA  )$ and $(\huaB, N, \rho_\huaB  )$ be two Lie
algebroids and let $\basemap: M\lon N$ be a smooth map.
\begin{itemize}

\item[1)] For a bundle map $\Psi:\basemap^! \huaB\lon \huaA$, let
its graph be
$$\Graph_{(\Psi,\phi)}=\set{(\Psi(v),v)|x\in M, v\in \huaB_{\basemap(x)}}
\subset \huaA\oplus \basemap^! \huaB.$$ Then,
$(\huaB,N)\RightLeft{\Psi}{\basemap}(\huaA,M)$ is a comorphism of
Lie algebroids if and only if $\Graph_{(\Psi,\phi)}$ is contained in
$\huaA\oplus_{\basemap}\huaB$ and $\Gamma(\Graph_{(\Psi,\phi)})$ is
a Lie subalgebra of $\Gamma(\huaA\oplus_{\basemap}\huaB)$.

\item[2)] For a bundle map $\Psi:\huaA\lon \huaB$, let its graph
be
$$\Graph_{(\Psi,\phi)}=\set{(u,\Psi(u))|x\in M, u\in \huaA_{x}}\subset
\huaA\oplus \basemap^! \huaB.$$ Then,
$(\huaA,M)\RightRight{\Psi}{\basemap}(\huaB,N)$ is a morphism of Lie
algebroids if and only if $\Graph_{(\Psi,\phi)}$ is contained in
$\huaA\oplus_{\basemap}\huaB$ and $\Gamma(\Graph_{(\Psi,\phi)})$ is
a Lie subalgebra of $\Gamma(\huaA\oplus_{\basemap}\huaB)$. \qed
\end{itemize}
\end{thm}
\begin{rmk}\label{Rmk:ItMayHappen}\rm
Evidently a comorphism of Lie algebroids
$(\huaB,N)\RightLeft{\Psi}{\basemap}(\huaA,M)$ requires that
$${\basemap_*}Im(\rho_\huaA|_{x})\supset Im(\rho_\huaB|_{\basemap(x)}).$$
And a morphism of Lie algebroids
$(\huaA,M)\RightRight{\Psi}{\basemap}(\huaB,N)$ requires that
$${\basemap_*}Im(\rho_\huaA|_{x})\subset Im(\rho_\huaB|_{\basemap(x)}).$$
\end{rmk}

In the following examples, we show that the traditional
representation and action theories can be restated in the language
of morphisms and comorphisms.
\begin{ex}\rm\label{Ex:derivativeRep}
Let $V\lon M$ be a vector bundle. Then one has the bundle of
\emph{covariant differential operators}, written
$\mathcal{CDO}(V)\lon M$ (c.f. \cite[III]{first},
\cite{Mackenzie:1995}, see also \cite{KSK}, where the notation
$\frkD(V)$ is used instead of $\mathcal{CDO}(V)$). A
\emph{derivative representation} of a Lie algebra $\LieG$ on $V$, is
a morphism of $\LieG$ into the Lie algebra
$\Gamma(\mathcal{CDO}(V))$ \cite{KSK}. This is obviously a
comorphism of Lie algebroids from $\LieG$ to $\mathcal{CDO}(V)$ ,
over the trivial map $M\lon pt$.
\end{ex}

\begin{ex}\rm\label{Ex:infitactionofLieaglebroid} \cite{KSK}
Let $(\huaA,M)$ be a Lie algebroid, and let $\varphi: Z\lon M$ be
a fibred manifold (i.e., $\varphi$ is a surjective submersion onto
$M$). An \emph{infinitesimal action of $\huaA$ on $Z$} is an
$\Real$-linear map $\Gamma(\huaA)\lon \mathcal{X}(Z)$, $A\mapsto
A_Z$, where $A\in\Gamma(\huaA)$, such that
\begin{itemize}
\item[1)] $A_Z$ is projective to $\rho_\huaA(A)$ (i.e., they are
$\varphi$-related); \item[2)] the map $A\mapsto A_Z$ preserves
brackets; \item[3)] the map $A\mapsto A_Z$ is $\CWM$-linear.
\end{itemize}

Of course this is equivalently saying that $(\cdot)\lon (\cdot)_Z$
over the fibred map $\varphi: Z\lon M$, is a comorphism of Lie
algebroids from $\huaA$ to $TZ$.
\end{ex}
\begin{ex}\rm\label{Ex:derivativeRepLAD} \cite{KSK} The derivative
representations of a Lie algebra in Example \ref{Ex:derivativeRep}
can be generalized to a Lie algebroid. Let $(\huaA,M)$ be a Lie
algebroid, and let $\varphi: Z\lon M$ be a fibred manifold and
$V\lon Z$ be a vector bundle on $Z$. A \emph{derivative
representation of $\huaA$ on $V$} associated to a given
infinitesimal action of $\huaA$ on $Z$: $A\mapsto A_Z$ (as in
Example \ref{Ex:infitactionofLieaglebroid}), is a morphism,
$\gamma$, of Lie algebras from $\Gamma(\huaA)$ to
$\Gamma(\mathcal{CDO}(V))$ such that
\begin{itemize}
\item[1)] for any $A\in \Gamma(\huaA)$,
$\rho_{{}_{\mathcal{CDO}(V)}}\gamma(A)=[\gamma(A),~\cdot~]=A_Z$;
\item[2)] $\gamma$ is $\CWM$-linear in the sense that
$\gamma(fA)=(\varphi^*f)\gamma(A)$, $\forall f\in\CWM$.
\end{itemize}

We point out that, such a derivative representation is a comorphism
of Lie algebroids
$$(\huaA,M)\RightLeft{\gamma}{\varphi}(\mathcal{CDO}(V),Z),\quad
\mbox{where we view  }{\gamma}\mbox{ as a map }:\ \varphi^!\huaA\lon
\mathcal{CDO}(V).
$$
Conversely, if a Lie algebroid comorphism is given as above, there
are associated
\begin{itemize}
\item[1)] an infinitesimal action of $\huaA$ on $Z$, $A\mapsto
[\gamma(A),~\cdot~]$; \item[2)] a derivative representation of
$\huaA$ on $V$, $A\mapsto \gamma(A)$, associated to the action of
 1).
\end{itemize}
\end{ex}
In fact, given any comorphism of Lie algebroids which is over a
fibred map, we can determine an action of Lie algebroids
\cite{PJHK2} (see also similar conclusions in Theorem
\ref{Thm:morphismGRDassAction}).
\begin{thm}\label{Thm:morphismLADassAction}
Let $\varphi: Z\lon M$ be a fibred manifold, and let $(\huaC,Z)$,
$(\huaA,M)$ be two Lie algebroids. If
$$(\huaA,M)\RightLeft{\Psi}{\varphi}(\huaC,Z)$$
is a comorphism of Lie algebroids, then the map
$$A\mapsto \rho_\huaC\circ\widetilde{\Psi}(A),
\quad \forall A\in \Gamma(\huaA),$$ defines an infinitesimal
action of $\huaA$ on $Z$. \qed
\end{thm}

\begin{ex}\label{Ex:ActionAlgebroid}\rm
[The action algebroid] Let $\LieG$ be a Lie algebra, $M$ a smooth
manifold and $\theta: \LieG\lon\XM$ a morphism of Lie algebras
(called an infinitesimal action of $\LieG$ on $M$). Then the vector
bundle $\huaA=M\times\LieG$ admits a Lie algebroid structure by
setting
\begin{eqnarray*}
[fX,gY] &=& fg[X,Y]+f\theta(X)(g)Y-g\theta(Y)(f)X,\\
\rho(fX)  & = & f\theta(X),
\end{eqnarray*}
$\forall f,g\in\CWM$, $X,Y\in\LieG$. We call $\huaA=M\times\LieG$
the action Lie algebroid induced by the action of $\LieG$ on $M$.
(This is a special case of Examples
\ref{Ex:infitactionofLieaglebroid} and
\ref{Ex:ActionAlgebroidExtended}.) Evidently, the action algebroid
$M\times\LieG$ admits both a trivial morphism to $\LieG$, and a
trivial comorphism from $\LieG$, both over $M\lon pt$.
\end{ex}

\begin{ex}\rm\label{Ex:ActionAlgebroidExtended}
We continue the assumptions in Example
\ref{Ex:infitactionofLieaglebroid}. It is shown in \cite{PJHK} that,
there is an associated Lie algebroid structure on $\varphi^!\huaA$
with the base $Z$, which is called the action Lie algebroid
associated to the action $A\lon A_Z$. And the diagram
\begin{equation}
\begin{CD}
\varphi^!\huaA @> \Id>>   \huaA\\
@V p VV @VV q V\\
Z @> \varphi>>   M
\end{CD}
\end{equation}
is a morphism of Lie algebroids. At the same time, there is a
comorphism $(\huaA,M)\RightLeft{\Id}{\varphi}(\varphi^!\huaA,Z)$.

It is also worth noting that  the action Lie algebroid
$(\varphi^!\huaA, Z)$ can be embedded into the $\varphi$-sum of
$(TZ,Z)$ and $(\huaA,M)$, by sending $(x,v)$ to $(v_Z,v)$, for each
$x\in Z$, $v\in \huaA|_{\varphi(x)}$.
\end{ex}

\section{Morphisms and Comorphisms of Lie Groupoids}
\label{Sec:GeometricModel}

This section is an exposition of the theory of two kinds of
morphisms concerning Lie groupoids, analogous to that of Lie
algebroids. We recall here some well-known definitions and certain
properties of Lie groupoids and their tangent Lie algebroids, and we
refer to \cite{AKA, first,Mackenzie:1995, Mkz:GTGA} for more details
(note that the composition convention in some of these texts is the
opposite of that followed here).

A groupoid $\Gamma$ on the base $M$, with respectively source and
target maps $\alpha$, $\beta$, will be denoted by $(\Gamma
\alphabetaarrow{M})$, or, more briefly, $(\Gamma, M)$.   We adopt
the convention that, whenever we write a product $gh$, we are
assuming that is defined, i.e., $\beta(g)=\alpha(h)$.   The base
$M\subset \Gamma$ will be regarded as the set of identities. The
inversion map $\iota: \Gamma\lon \Gamma$ will be denoted by
$\iota(g)=g\inverse$.

For $x\in M$, its \textbf{orbit}, denoted by $O_x$, is the set
$\beta\circ\alpha\inverse(x)\subset M$.

Let $(\Pi,N)$ be a groupoid and $\basemap: M\lon N$ be a map. We
define the $\basemap$-pullback of $\Pi$ (with respect to the
$\alpha$-fiber)
$$
M\times_{\basemap}\Pi\defbe \set{(x,{w})| x\in M, \basemap(x)=
\alpha_\Pi (w)},
$$
which is also denoted by $\basemap^!\Pi$, or $\basemap^*\Pi$ by many
authors. When $\Pi$ is a Lie groupoid on $N$ and $\basemap$ is a
smooth map between smooth manifolds, $\alpha_\Pi$ is a submersion
and hence $M\times_{\basemap}\Pi$ is a smooth manifold, and $\pr_M$,
the projection to $M$ is a submersion.

We are now ready to introduce the two concepts of morphism of
groupoids analogous to   the morphisms and comorphisms of Lie
algebroids. We will prove that they are global formulations in terms
of Lie groupoids. In \cite{PJHK2}, Higgins and Mackenzie had given
the definition of comorphisms of groupoids in the language of
actions. Here we prefer to adopt a direct description as follows.

\begin{defi}\label{Def:ComorphismGroupoid}
Let $(\Gamma \rightrightarrows{M},\alpha_\Gamma, \beta_\Gamma)$ and
$(\Pi \rightrightarrows{N},\alpha_\Pi,\beta_\Pi)$ be two groupoids
on bases $M$ and $N$ respectively. A \textbf{comorphism of
groupoids} from $\Pi$ to $\Gamma$, over a given map $\basemap: M\lon
N$, is a map $\Phi:M\times_{\basemap}\Pi\lon \Gamma$, written
$(\Pi,N)\RightLeft{\Phi}{\basemap}(\Gamma,M)$, such that the diagram
\begin{equation}\label{Graph:phizeroPullback}
\begin{CD}
~~M\times_{\basemap}\Pi @> {\Phi}>>   ~~~~\ \ \Gamma\\
{\pr_M}@VV V {\alpha_\Gamma} @VV
 V\\
~~M @> \Id_M>>   ~~~~\ \ M
\end{CD}
\end{equation}
commutes and the following conditions hold
\begin{itemize}
\item[1)] for all $x\in M$, $\Phi(x,\basemap(x))=x$;

\item[2)] for all $(x,{w})\in M\times_{\basemap}\Pi$,
$\basemap\circ\beta_\Gamma\circ\Phi(x,{w})=\beta_\Pi({w})$;

\item[3)] for all $(x,{w})\in M\times_{\basemap}\Pi$,
$(\beta_{\Gamma}\circ\Phi(x,{w}),{z})\in M\times_{\basemap}\Pi$,
there holds
$${\Phi}(x,{w}{z}) =
\Phi(x,{w})\Phi(\beta_{\Gamma}\circ\Phi(x,{w}),{z}).$$
\end{itemize}

In particular, if $\basemap$ is surjective and $\Phi$ is  injective,
then we call $(\Pi,N)$ a \textbf{co-subgroupoid} of $(\Gamma,M)$.
\end{defi}
The morphisms of groupoids are already a well-known concept
\cite{first}, and they are global version of the morphisms of Lie
algebroids.
\begin{defi}\label{Def:MorphismGroupoid}
Let $(\Gamma \rightrightarrows{M},\alpha_\Gamma, \beta_\Gamma)$
and $(\Pi \rightrightarrows{N},\alpha_\Pi,\beta_\Pi)$ be two
groupoids on bases $M$ and $N$ respectively. A \textbf{morphism of
groupoids} over $\basemap: M\lon N$ is a map $\Phi$: $\Gamma\lon
\Pi$, written $(\Gamma,M)\RightRight{\Phi}{\basemap}(\Pi,N)$, such
that
\begin{itemize}
\item[1)] $\Phi(M)\subset N$ and $\basemap=\Phi|_M$; \item[2)]
$\alpha_\Pi \circ\Phi=\basemap\circ\alpha_\Gamma$,
$\beta_\Pi\circ\Phi=\basemap\circ\beta_\Gamma$; \item[3)]
$\Phi(gh)=\Phi(g)\Phi(h)$, for all composable $g, h\in \Gamma$.
\end{itemize}

In particular, if $\basemap$   and $\Phi$ are both injective, then
we call $(\Gamma,M)$ a \textbf{subgroupoid} of $(\Pi,N)$.
\end{defi}

For the base-preserving case
$M\stackrel{\basemap=\Id}{\longrightarrow} N=M$, these two kinds of
morphisms coincide.

A comorphism of \textsl{Lie groupoids}
$(\Pi,N)\RightLeft{\Phi}{\basemap}(\Gamma,M)$, or a morphism
$(\Gamma,M)\RightRight{\Phi}{\basemap}(\Pi,N)$, is defined
similarly, with the additional requirement that ${\Phi}$, as well as
$\basemap$ in the above definitions to be smooth maps.

For a Lie groupoid $(\Gamma\alphabetaarrow{M})$, we define the
\emph{tangent Lie algebroid} $(Lie\Gamma,\rho)$ as
$$Lie\Gamma\defbe \bigcup_{x\in M}T_x\alpha\inverse(x)=
\set{v\in T_x\Gamma| x\in M, \alpha_*(v)=0}.$$ The bracket of
$\Gamma(Lie\Gamma)$ is determined by the commutator of
left-invariant vector fields and the anchor map is given by
$\rho=\beta_{\Gamma*}|_M$.

The infinitesimal counterpart of a comorphism can be determined by
differentiation, yielding a comorphism of Lie algebroids. This fact
is illustrated by the following theorem.

\begin{thm}Let $(\Pi,N)\RightLeft{\Phi}{\basemap}(\Gamma,M)$
be a comorphism of Lie groupoids. Then, the tangent map of $\Phi$
induces a vector bundle morphism $\Psi=\Phi_*$
$:\basemap^!(Lie\Pi)\lon Lie\Gamma$ and $(Lie
\Pi,N)\RightLeft{{\Psi}}{{\basemap}}(Lie \Gamma,M)$ is a
comorphism of Lie algebroids.
\end{thm}
\pf By  (1) of Definition \ref{Def:ComorphismGroupoid}, $\Phi_*$
sends $T_{(x,\basemap(x))}\alpha_\Pi\inverse(\basemap(x))$ to
$T_{x}\alpha_\Gamma\inverse(x)$ and hence $\Psi|_{x}\defbe
{\Phi_*}|_{(x,\basemap(x))}$, for all $x\in M$, defines a bundle
map. According to (2), we have
$${\basemap_*}\circ\beta_{\Gamma*}\circ\Phi_*(x,v)=\beta_{\Pi*}(v),
\quad\forall v\in Lie\Pi_{\basemap(x)},$$ i.e., $\Psi$ is subject to
the condition ${\basemap_*}\circ\rho_{Lie\Pi}\circ\Psi=
\rho_{Lie\Gamma}$. It remains to prove that the induced map
$\widetilde{\Psi}: \Gamma(Lie\Pi)\lon \Gamma(Lie\Gamma)$ is a
morphism of Lie algebras.

For each section $B\in \Gamma(Lie\Pi)$, we denote the corresponding
left-invariant vector field on $\Pi$ by $\LeftMove{B}$ and
similarly, $\LeftMove{\widetilde\Psi(B)}$ denotes the left-invariant
vector field on $\Gamma$. Let us prove that $\LeftMove{B}$, regarded
as a vector field on $M\times_{\basemap}\Pi$,  is $\Phi$-related to
$\LeftMove{\widetilde\Psi(B)}$.   For each $x\in M$,
${w}\in\alpha_\Pi\inverse(\basemap(x))$, if we suppose that
$B_{\beta_\Pi({w})}=\frac{d}{dt}c(0)$, where $c: [0,1]\lon
\alpha_\Pi\inverse(\beta_\Pi({w}))$ is a smooth curve, then
\begin{eqnarray*}
&&\Phi_{*(x,{w})}\LeftMove{B }({w})\\
&=& \frac{d}{dt}|_{t=0}\Phi(x,{w}c(t)) \\&=&
\frac{d}{dt}|_{t=0}\Phi(x,{w})\Phi(\beta_\Gamma\circ\Phi(x,{w}),c(t))
\quad\mbox{(by (3) in Definition \ref{Def:ComorphismGroupoid})}\\
&=&L_{\Phi(x,{w})*}(\frac{d}{dt}|_{t=0}\Phi(\beta_\Gamma\circ\Phi(x,{w}),c(t))),
\\
&=&L_{\Phi(x,{w})*}(\Psi(\beta_\Gamma\circ\Phi(x,{w}),B_{\beta_\Pi({w})}))=
\LeftMove{\widetilde\Psi(B)}(\Phi(x,{w})).
\end{eqnarray*}
Hence for two $B_1$, $B_2\in \Gamma(Lie\Pi)$,
$[\LeftMove{B_1},\LeftMove{B_2}]=\LeftMove{[B_1,B_2]}$ is also
$\Phi$-related to
$[\LeftMove{\widetilde\Psi(B_1)},\LeftMove{\widetilde\Psi(B_2)}]$
$=$  $\LeftMove{[\widetilde\Psi(B_1),\widetilde\Psi(B_2)]}$. It
follows that
$$\widetilde\Psi([B_1,B_2])=
\Phi_*|_M(\LeftMove{[B_1,B_2]})=
\LeftMove{[\widetilde\Psi(B_1),\widetilde\Psi(B_2)]}|_M=
[\widetilde\Psi(B_1),\widetilde\Psi(B_2)].
\quad\quad\blacksquare$$

Similarly, for a morphism
$(\Gamma,M)\RightRight{\Phi}{\basemap}(\Pi,N)$ of Lie groupoids,
the tangent map $\Phi_*$ which evidently sends $Lie \Gamma$ to
$Lie \Pi$, is a morphism of Lie algebroids (over $\basemap$). The
details of the proof can be found in \cite{PJHK}.

\noindent\textbf{$\bullet$ The $\basemap$-product of groupoids}.\\
There is likewise a global version of the restriction theory for
groupoids. Let $(\Gamma \alphabetaarrow{M})$ be a groupoid.
$M_0\subset M$ a subset. We call
$$\Gamma_{M_0}\defbe\set{g\in \Gamma| \alpha(g)\in M_0,\beta(g)\in M_0}$$
the $M_0$-\textbf{restriction} of $\Gamma$. It is easy to see that
$\Gamma_{M_0}$ is also a groupoid\footnote{For a Lie groupoid
$\Gamma$, the restriction $\Gamma_{M_0}$ is not necessarily a Lie
groupoid.} on the base $M_0$, inheriting the structures of $\Gamma$
(c.f. Example \ref{Ex:QuotientReduction}).

For two groupoids $(\Gamma \rightrightarrows{M},\alpha_\Gamma,
\beta_\Gamma)$ and $(\Pi \rightrightarrows{N},\alpha_\Pi,\beta_\Pi)$
on bases $M$ and $N$ respectively, we can endow their direct product
$\Gamma\times \Pi$ with a groupoid structure on $M\times N$ by
setting
$$\alpha(g,{w})\defbe (\alpha_\Gamma(g),\alpha_\Pi({w})),\quad
\beta(g,{w})\defbe (\beta_\Gamma(g),\beta_\Pi({w}));$$
$$(g,{w})(h,{z})\defbe (gh,{w}{z}),\quad \mbox{where } g,h\in\Gamma \mbox{ and } {w},{z}\in \Pi
\mbox{ are respectively composable}.$$ We will denote this groupoid
by $\Gamma\widetilde{\times} \Pi$, and call it the \textbf{direct
product} of $\Gamma$ and $\Pi$.


\begin{defi}\label{Def:phi0sum}
The $\Graph(\basemap)$-restriction of the direct sum groupoid
$\Gamma\widetilde{\times}\Pi$, denoted by
$\Gamma\widetilde{\times}_{\basemap}\Pi$, is called the
\textbf{$\basemap$-product} of $\Gamma$ and $\Pi$.
\end{defi}
One may characterize such a groupoid by
$$\Gamma\widetilde{\times}_{\basemap}\Pi\defbe
\set{(g,{w})\in \Gamma\widetilde{\times}\Pi | \alpha_\Pi({w}) =
\basemap\circ\alpha_\Gamma(g), \beta_\Pi({w}) =
\basemap\circ\beta_\Gamma(g)}.
$$

Notice that, $\Gamma\widetilde{\times}_{\basemap}\Pi$ is a
groupoid\footnote{The $\basemap_{0}$-product of two Lie groupoids is
not necessarily a Lie groupoid.} on base space $M$. The source and
target maps send $(g,{w})\in \Gamma\widetilde{\times}_{\basemap}\Pi$
to $\alpha_\Gamma(g)$ and $\beta_\Gamma(g)$ respectively (c.f
Example \ref{Ex:DirectSum}).

On the global level, the analogues of Theorem
\ref{Thm:morphismLPA} and Theorem \ref{Thm:coandmorphismLAB} are
as follows. Again we omit the proof.
\begin{thm}\label{Thm:coandmorphismLGB} Let $(\Gamma \rightrightarrows{M},\alpha_\Gamma, \beta_\Gamma)$
and $(\Pi \rightrightarrows{N},\alpha_\Pi,\beta_\Pi)$ be two
groupoids. Let $\basemap: M\lon N$ be a map.
\begin{itemize}
\item[1)] A map $\Phi: M\times_{\basemap}\Pi\lon\Gamma$
is a comorphism of groupoids if and only if the graph
$$\Graph_{(\Phi,\phi)}\defbe \set{(\Phi(x,g),g)|x\in M,
g\in \alpha_\Pi\inverse(\basemap(x))}
$$ is a subgroupoid of $\Gamma\widetilde{\times}_{\basemap}\Pi$,
on the same base $M$.

\item[2)] A map $\Phi: \Gamma\lon \Pi$
is a morphism of groupoids if and only if the graph
$$\Graph_{(\Phi,\phi)}\defbe \set{(g,\Phi(g))|g\in \Gamma}
$$ is a subgroupoid of $\Gamma\widetilde{\times}_{\basemap}\Pi$,
on the same base $M$. \qed
\end{itemize}
\end{thm}
\begin{rmk}\rm
Evidently a comorphism of groupoids
$(\Pi,N)\RightLeft{\Phi}{\basemap}(\Gamma,M)$ requires that the
$\basemap$-image of each orbit of $\Gamma$ covers some orbit of
$\Pi$, i.e., $O_{\basemap(x)}\subset \basemap( O_{x})$. And a
morphism of groupoids $(\Gamma,M)\RightRight{\Phi}{\basemap}(\Pi,N)$
requires that $\basemap$ sends each orbit of $\Gamma$ into an orbit
of $\Pi$, i.e., $\basemap( O_{x})\subset O_{\basemap(x)}$ (c.f.
Remark \ref{Rmk:ItMayHappen} ).
\end{rmk}

\noindent\textbf{$\bullet$  Examples.}\\
In the remaining part of the paper, we relate some examples showing
various kinds of morphisms and comorphisms of Lie algebroids and
groupoids, to be compared with the examples of actions and
representations in the preceding section.

\begin{ex}\label{Ex:pairgroupoid}\rm Let $M$ be a smooth manifold.
Then $M\times M$ admits a groupoid structure, called the
\emph{pair groupoid} on $M$. The tangent Lie algebroid for
$M\times M$ is $TM$ \cite{AKA}.
\end{ex}
\begin{ex}\label{Ex:tangentComorphism}\rm Let $\basemap: M\lon N$ be a smooth map
between smooth manifolds. The map $\basemap: M\lon N$ naturally
induces  a morphism of the pair groupoids
$$(M\times M,M)\RightRight{\Phi}{\basemap}(N\times N,N),$$
by setting $\Phi(x,y)=(\basemap(x),\basemap(y))$, $\forall x,y\in
M$. Clearly, ${\basemap_*}: TM\lon TN$ is the corresponding morphism
of their tangent  Lie algebroids.
\end{ex}
\begin{ex}\rm We investigate the comorphisms of the pair groupoids.
Suppose that a map  $\basemap: M\lon N$ is given and there is a
comorphism $\Phi: M\times_{\basemap} (N\times N)\lon M\times M$ over
$\basemap$. Here we see that
$$
M\times_{\basemap} (N\times N)= \set{(x,\basemap(x),y)| x\in M,y\in
N}.
$$

According to Definition \ref{Def:ComorphismGroupoid}, $\Phi$ can be
written in  the form
$$
\Phi: (x,\basemap(x),y)\mapsto (x,\Phi_0(x,y)),\quad\forall x\in
M,y\in N,
$$
where $\Phi_0: M\times N\lon M$. And the axioms in Definition
\ref{Def:ComorphismGroupoid} are now expressed as follows:
\begin{itemize}
\item[1)] $\Phi_0(x,\basemap(x))=x$, $\forall x\in M$;
\item[2)] $\basemap\circ\Phi_0(x,y)=y$, $\forall x\in M$, $y\in N$;
\item[3)] $\Phi_0(x,z)=\Phi_0(\Phi_0(x,y),z)$, $\forall x\in M$, $y,z\in N$.
\end{itemize}
Thus $\basemap$ must be a surjection. And for each fixed $x\in M$,
$\Phi_0(x,~\cdot~): N\lon M$ is an injection as well as a right
inverse of $\basemap$. If $M$ and $N$ are both smooth manifolds,
these conditions are exactly saying that there is a foliation
structure on $M$ and each leaf is diffeomorphic to $N$ via
$\basemap$. So generally speaking, it is hard to find a comorphism
from $\Pi=N\times N$ to $\Gamma=M\times M$.
\end{ex}
\begin{ex}\rm
We now study what a comorphism  $\Psi$, from the algebroid $TN$ to
$TM$, could be. Suppose that $\Psi$ is over $\basemap$, then by
Remark \ref{Rmk:ItMayHappen}, $Im(\basemap)$ must be an open
submanifold of $N$, and $\basemap$ must be a submersion from $M$ to
$Im(\basemap)$. And the map $\widetilde{\Psi}$ must assign every
$Y\in \XN$ to a lifted $\widetilde{\Psi}(Y)\in \XM$, such that
${\basemap_*}(\widetilde{\Psi}(Y))=Y$ and the bracket must be
preserved
$$\widetilde{\Psi}{[X,Y]}=[\widetilde{\Psi}(X),\widetilde{\Psi}(Y)].$$
For each $x\in M$, let $\huaF_x=Im(\Psi(\set{x}\times
T_{\basemap(x)}N))$. The conditions above imply that   each
$\huaF_x$ is isomorphic to $T_{\basemap(x)}N$ via $\basemap_*$. The
distribution $\huaF$ on $M$ is integrable. And each leaf of $\huaF$
is a covering space of $Im(\basemap)$.
\end{ex}
\begin{ex}\rm Let $M$ be a submanifold of $N$ and $\huaT\subset
TM$ a regular, integrable distribution. Then the inclusion
$\huaT\subset TN$ together with the embedding $M\subset N$ is an
injective comorphism of Lie algebroids. On the other hand, we have a
groupoid with base $M$,
$$M\times_\huaT M\defbe \set{
(x,y)\in M\times M| x \mbox{ and } y \mbox{ belong to the same
integral submanifold of }M,\mbox{ via }\huaT }.$$ The embedding of
$M\times_\huaT M$ into $N\times N$ is a morphism of groupoids.
\end{ex}

\begin{ex}\label{Ex:ActionGroupoid}\rm[The action groupoid]
Let $G$ be a group with the unit element $e$ and let $G$ act on a
set $M$ (to the right), i.e., we have a map $\Theta: M\times G\lon
M$, $(x,g)\mapsto xg$, satisfying the axioms: $xe=x$,
$x(g_1g_2)=(xg_1)g_2$, for all $x\in M$, $g_1$, $g_2\in G$. Then
$\Gamma=M\times G$ admits a groupoid structure with base space $M$,
as follows
\begin{eqnarray*}
\alpha: \quad (x,g)&\mapsto& x, \\
\beta:  \quad (x,g)&\mapsto& xg,\\
\iota:  \quad (x,g)&\mapsto& (xg,g\inverse), \\
(x,g_1)(xg_1,g_2)&=&(x,g_1g_2).
\end{eqnarray*}
We call $\Gamma=M\times G$ the action groupoid induced by the action
of $G$ on $M$. Then $Lie\Gamma$ is exactly the action algebroid
$M\times Lie G$ given in Example \ref{Ex:ActionAlgebroid}, where
$Lie G=T_eG$ has the Lie bracket coming from left-invariant vector
fields. This is a particular case of Example
\ref{Ex:ActionGroupoidExtended}, \ref{Ex:ActionGroupoidEX}.
\end{ex}

\begin{ex}\label{Ex:semilinearRepLieGroup}
\rm Given a vector bundle $V\lon M$, we denote by $\mathcal{LF}(V)$
the \textsl{linear frame groupoid} of $V$ (see \cite[III]{first},
where it is denoted by $\Pi(V)$), which is the collection of all
linear isomorphisms from a fiber of $V$ to some generally different
fiber of $V$. The Lie algebroid of $\mathcal{LF}(V)$ is canonically
isomorphic to $\mathcal{CDO}(V)$ (this fact was discovered
independently by Kumpera \cite{Kumpera}, Mackenzie \cite{first},
Hermann \cite{Hermann}, Kosmann-Schwarzbach \cite{K1976}). The
 group of bisections $\mathcal{LF}(V)$ is in fact $Aut(V)$, the group
of vector bundle automorphisms of $V$. (One can also regard it as
the \emph{semi-linear isomorphisms} of $\Gamma(V)$, see
\cite{K1976,KSK}.) In \cite{KSK}, Kosmann-Schwarzbach and  Mackenzie
defined the \emph{semi-linear representation} of a Lie group $G$ on
$V$ to be a group morphism, say $R$, from $G$ to $Aut(V)$, with the
smooth condition that, the map $(x,g)\mapsto R(g)(\nu)(x)$, $M\times
G\lon V$, is smooth, for each $\nu\in\Gamma(V)$ (or equivalently,
the map $(x,g)\lon R(g)|_{x}$, $M\times G\lon \mathcal{LF}(V)$ is
smooth).

This is actually a comorphism of Lie groupoids from $G$ to
$\mathcal{LF}(V)$ (over $M\lon pt$), and it provides a global
formulation for derivative representations of Lie algebras
(Example \ref{Ex:derivativeRep}).
\end{ex}
\begin{ex}\rm\label{Ex:ActionGroupoidExtended}
\cite{PJHK,KSK} Let $(\Gamma,M)$ be a groupoid, and let $\varphi:
Z\lon M$ be a surjective map. A (\emph{right}) action of $\Gamma$
on $Z$ is a map $S$ assigning each $g$ of $\Gamma$ to another map
$$S(g): Z_{\alpha(g)}\lon Z_{\beta(g)},$$
such that
\begin{itemize}
\item[1)] $S(y)=\Id: Z_y\lon Z_y$, for any $y\in M$. (It seems that
 \cite{KSK}  has missed this condition.) \item[2)] $S(g_1g_2)=S(g_2)\circ S(g_1)$,
for any composable $g_1$, $g_2\in \Gamma$.
\end{itemize}
We argue that, an action of $\Gamma$ is   a comorphism of groupoids
from $\Gamma$ to the pair groupoid $Z\times Z$. In fact, if the
action $S$ is given as above, then by setting
$$\Phi:{Z\times _{\varphi}\Gamma}\lon Z\times Z,\quad
(x,g)\mapsto (x,S(g)(x)),\quad\forall x\in Z,g\in
\alpha_\Gamma\inverse(\varphi(x)),
$$
we obtain a comorphism $\Phi$. Conversely, given $\Phi$, we can
determine the action $S$ by
$$S(g)(x)\defbe \pr_2\circ\Phi(x,g), \quad\mbox{where }\pr_2
 \mbox{ is the projection to the second component.}$$

When considering a Lie groupoid $\Gamma$ and a fibred manifold
$\varphi: Z\lon M$, we require such an action to be a smooth
comorphism of Lie groupoids. At the same time, one obtains an
infinitesimal action of the corresponding Lie algebroid
$Lie\Gamma$ on $Z$. For more details about these groupoids and
their infinitesimal forms, see   Examples \ref{Ex:pairgroupoid},
\ref{Ex:ActionAlgebroid}, \ref{Ex:ActionGroupoid} and
\ref{Ex:ActionAlgebroidExtended}.
\end{ex}
The theorem which follows is an analogue of Theorem
\ref{Thm:morphismLADassAction}: for each comorphism of Lie
groupoids over a fibred map, there is associated an action of a
Lie groupoid.
\begin{thm}\label{Thm:morphismGRDassAction}
Let $\varphi: Z\lon M$ be a fibred manifold, and let $(\Omega,Z)$,
$(\Gamma,M)$ be two Lie groupoids. If
$$(\Gamma,M)\RightLeft{\Phi}{\varphi}(\Omega,Z)$$
is a comorphism of Lie groupoids, then the map $g\mapsto S(g)$,
where
$$S(g): Z_{\alpha(g)}\lon Z_{\beta(g)},\; x\mapsto
\beta_\Omega\circ\Phi(x,g), \quad \forall x\in Z_{\alpha(g)},$$
defines a Lie groupoid action of $\Gamma$ on $Z$. \qed
\end{thm}

\begin{ex}\rm\label{Ex:ActionGroupoidEX}
Recall Example \ref{Ex:ActionGroupoidExtended}. It was shown in
\cite{Ehresmann} (see also \cite{Cranic, KSK}) that whenever an
action of a groupoid $(\Gamma,M)$ on $\varphi: Z\lon M$ is defined,
there is an associated groupoid structure on the pull-back space
$Z\times_\varphi\Gamma$, with base $Z$, known as the \emph{action
groupoid} associated to the given action. The source map is
$(x,g)\mapsto x$, while the target map is $(x,g)\mapsto S(g)(x)$.
The multiplication is $(x,g)(S(g)(x),h)=(x,gh)$. There is obviously
a morphism from $Z\times_\varphi\Gamma$ to $\Gamma$, $(x,g)\mapsto
g$, over $\varphi$.
\end{ex}

\begin{ex}\rm
In Example \ref{Ex:derivativeRepLAD}, we have discussed the
derivative representations of Lie algebroids. On the global level,
a \emph{semi-linear representation of a Lie groupoid} associated
to an action of $(\Gamma,M)$ on the fibred manifold
$Z\stackrel{\varphi}{\lon} M$ (defined as in Example
\ref{Ex:ActionGroupoidExtended}), is a morphism $\Phi$ from the
action Lie groupoid $Z\times_\varphi\Gamma$ to $\mathcal{LF}(V)$
(on the same base $Z$), where $V\lon Z$ is a vector bundle. This
notion was studied in \cite{Brown}, and also presented in
\cite{KSK}. Here we emphasize that  such a semi-linear
representation is in fact a comorphism of Lie groupoids,
$$(\Gamma,M)\RightLeft{\Phi}{\varphi}(\mathcal{LF}(V),Z).$$
And the associated action of $\Gamma$ on $Z$ is determined by the
relation
$$S(g)(x)=\beta_{\mathcal{LF(V)}}\circ\Phi(x,g),\quad
\forall x\in Z,g\in\alpha_\Gamma\inverse(\varphi(x)).$$
\end{ex}

\begin{ex}\rm
Recall the action groupoid $Z\times_\varphi\Gamma$ given in Example
\ref{Ex:ActionGroupoidEX}. We point out that  it is in fact a
subgroupoid of the $\varphi$-product of the pair groupoid $Z\times
Z$ and $\Gamma$. In fact, every element $(x,g)\in{Z\times
_{\varphi}\Gamma}$ can be regarded as $((x,S(g)x);g)\in (Z \times
Z)\widetilde{\times}_{\varphi}\Gamma$, where $x\in Z$, $g\in
\alpha_\Gamma\inverse(\varphi(x))$.
\end{ex}

\begin{ex}\label{Ex:ComorphismActiongroupoid}\rm
Let $G$, $H$ be two groups which act on $M$ and $N$ respectively.
Suppose that there is a map $\basemap: M\lon N$ and a group morphism
$\tau: G\lon H$, and that they are compatible in the following sense
\begin{equation}\label{Eqt:CoMorphismAction}
\basemap(xg)=\basemap(x)\tau(g),\quad \forall x\in M, g\in G.
\end{equation}
Then, the map between action groupoids $\Phi:M\times G\lon N\times
H$ given by
$$\Phi(x,g)=(\basemap(x),\tau(g)),\quad \forall x\in M, g\in G,$$
is a morphism of groupoids (over $\basemap$).
\end{ex}
\begin{ex}\label{Ex:MorphismActiongroupoid}\rm
Let $G$, $H$ be two groups which act respectively on $M$ and $N$.
Suppose that there is a map $\basemap: M\lon N$ and a group morphism
$\varsigma: H\lon G$, and that they are compatible in the following
sense
\begin{equation}\label{Eqt:MorphismAction}
\basemap(x)h=\basemap(x\varsigma(h)),\quad \forall x\in M, h\in H.
\end{equation}
Then, there is a comorphism of groupoids $(N\times
H,N)\RightLeft{\Phi}{\basemap}(M\times G,M)$, where $\Phi$ is
defined by
$$(x,(\basemap(x),h))\stackrel{\Phi}{\longrightarrow}
(x,\varsigma(h)),\quad\forall x\in M,h\in H.$$
\end{ex}

Note that the $\basemap$-product of the two action groupoids
$\Gamma=M\times G$ and $\Pi=N\times H$ in the preceding two examples
is
$$
\Gamma\widetilde{\times}_{\basemap}\Pi=\set{(x,g;\basemap(x),h)|x\in
M,g\in G,h\in H,\basemap(xg)=\basemap(x)h}.
$$

\begin{ex}\label{Ex:ComorphismActionalgebroid}\rm
Let $\LieG$ and $\LieH$ be two Lie algebras, and $M$, $N$ two smooth
manifolds. Let $\theta: \LieG\lon \XM$ and $\eta: \LieH\lon \XN$ be
two Lie algebra morphisms. Suppose that $\Upsilon: \LieG\lon \LieH$
is also a Lie algebra morphism, $\basemap: M\lon N$ a smooth map and
that they are compatible in the sense that
\begin{equation}\label{Eqt:CoMorphismActionalgebroid}
\eta(\Upsilon(X))|_{\basemap(x)}={\basemap_*}(\theta(X)|_{x}),\quad
\forall x\in M,X\in\LieG,
\end{equation}
(i.e., $\theta(X)$ and $\eta(\Upsilon (X))$ are $\basemap$-related.)
Then, we find a morphism of the action Lie algebroids $\Psi: M\times
\LieG\lon N\times \LieH$ (over $\basemap$), given by
$$\Psi(x,X)=(\basemap(x),\Upsilon(X)),\quad \forall x\in M,X\in\LieG.$$
\end{ex}
\begin{ex}\label{Ex:MorphismActionalgebroid}\rm
Let $\LieG$ and $\LieH$ be two Lie algebras, and $M$, $N$ two smooth
manifolds. Let $\theta: \LieG\lon \XM$, $\eta: \LieH\lon \XN$ be two
Lie algebra morphisms. Suppose that $\Sigma: \LieH\lon \LieG$ is
also a Lie algebra morphism, $\basemap: M\lon N$ a smooth map and
that they are compatible in the sense that
\begin{equation}\label{Eqt:MorphismActionalgebroid}
\eta(Y)|_{\basemap(x)}={\basemap_*}(\theta(\Sigma(Y))|_{x}),\quad
\forall x\in M,Y\in\LieH,
\end{equation}
(i.e., $\theta(\Sigma(Y))$ and $\eta(Y)$ are $\basemap$-related.)
Then, there is a comorphism of the action Lie algebroids $\Psi:
\Gamma(N\times \LieH)\lon \Gamma(M\times \LieG)$ (over $\basemap$),
given by
$$\Psi(fY)=\basemap^*(f)\Sigma(Y),\quad \forall f\in\CWN,Y\in\LieH.$$
\end{ex}
These two morphisms are the infinitesimal versions of Example
\ref{Ex:ComorphismActiongroupoid} and
\ref{Ex:MorphismActiongroupoid}.
\begin{ex}\label{Ex:GaugegroupAndalgebroid}\rm[The gauge groupoid]
 Let $(P\stackrel{p}{\lon}M, G)$
be a principal bundle on $M$ with  structure group $G$. Then $G$
acts diagonally on $P\times P$. We denote the orbit of $( x_1, x_2)$
by $\langle  x_1, x_2\rangle $ and the orbit manifold by ${(P\times
P)/G}$. Define $\alpha: {(P\times P)/G}\lon M$ by $\langle  x_1,
x_2\rangle \mapsto p( x_1)$ and similarly $\beta(\langle  x_1,
x_2\rangle )=p( x_2)$. For two elements $\langle  x_1, x_2\rangle $,
$\langle  y_1, y_2\rangle \in {(P\times P)/G}$, if $p( x_2)=p(
y_1)$, then there is a unique element $g\in G$ such that $ x_2= y_1
g$ and we define $\langle
 x_1, x_2\rangle \langle
 y_1, y_2\rangle =\langle  x_1, x_2\rangle \langle  y_1 g,  y_2 g\rangle
=\langle  x_1, y_2g\rangle $. In this way, ${(P\times P)/G}$ carries
a Lie groupoid structure with base $M$, called the gauge groupoid of
$(P,M,G)$ \cite{Mackenzie:1995,WeinSymplecticGrod}.

The action of $G$ on $P$ naturally lifts to an action on $TP$. We
denote the orbit of $X_x\in T_xP$ by $\langle X_x\rangle $ and the
quotient manifold by $\frac{TP}{G}$. Since this action is free,
$\frac{TP}{G}$ admits a vector bundle structure with base $M$, and
bundle projection $q: \langle X_x\rangle \mapsto p(x)$. Sections of
$\frac{TP}{G}$ can be regarded as vector fields on $P$ which are
$G$-invariant. It follows that $\Gamma(\frac{TP}{G})$ has an induced
bracket transferred from $\mathcal{X}({P})$. Besides, the tangent
map $p_*$ can also be transferred to $\frac{TP}{G}\lon TM$. Thus,
$(\frac{TP}{G},p_*,M)$ is a Lie algebroid, which is in fact the
tangent Lie algebroid for the gauge groupoid ${(P\times P)/G}$
\cite{Aty:Complexfiber,Mackenzie:1995}.
\end{ex}
\begin{ex}\label{Ex:morphismGaugealgebroid}\rm
Continuing the above example, the map $\langle ~\cdot~\rangle $
together with $p$ is obviously a morphism of Lie algebroids from
$(TP,\Id,P)$ to $(\frac{TP}{G},p_*,M)$. Like-wise, we get a typical
example of comorphism of Lie algebroids from $(\frac{TP}{G},p_*,M)$
to $(TP,\Id,P)$, by sending $X\in \Gamma(\frac{TP}{G})$ identically
to $X\in \mathcal{X}(P)$ which is $G$-invariant, over the projection
$p: P\lon M$.
\end{ex}
\begin{ex}\rm  Let $(P\stackrel{p}{\lon}M, G)$ and
$(Q\stackrel{q}{\lon}N, H)$ be two principal bundles on base spaces
$M$, $N$ respectively. Suppose that $\tau: G\lon H$ is a morphism of
Lie groups and $\basemap: M\lon N$ a smooth map. Given a smooth map
$\Phi: P\lon Q$ compatible with $\tau$ and $\basemap$ in the sense
that
$$q\circ \Phi(x)=\basemap\circ p(x),\quad
\Phi(xg)=\Phi(x)\tau(g),\quad \forall x\in M,g\in G,$$ we can then
construct a morphism of Lie groupoids $({(P\times
P)/G},M)\RightRight{\Phi}{\basemap}( {(Q\times Q)}/{H},N)$, given by
$$\Phi(\langle  x_1, x_2\rangle )=\langle \Phi( x_1),\Phi( x_2)\rangle ,\quad \forall  x_1, x_2\in P.$$
At the same time, we obtain a  morphism of Lie algebroids $\Phi_*$
from $(\frac{TP}{G},p_*,M)$ to $(\frac{TQ}{H},q_*,N)$, over
$\basemap$.
\end{ex}
\begin{ex}\rm  Let $(P\stackrel{p}{\lon}M, G)$ be a principal
bundle. Consider the pair groupoid $P\times P$ and the gauge
groupoid ${(P\times P)/G}$ given in Example \ref{Ex:pairgroupoid},
\ref{Ex:GaugegroupAndalgebroid}. There are canonically defined
morphism and comorphism of groupoids (both over $p$),
$$(P\times P,P)\RightRight{\langle \cdot,\cdot\rangle }{p}
(\frac{P\times P}{G},M), \quad \quad ({\frac{P\times
P}{G}},M)\RightLeft{\Phi}{p}(P\times P,P).$$ We elaborate on the
second one.  For any element in $p^!({(P\times P)/G})$, say
$(x,\langle y,z\rangle )$, where $x,y,z\in P$ and $p(x)=p(y)$, we
find $g\in G$ such that $y=xg$, and hence $(x,\langle y,z\rangle
)=(x,\langle x,zg\inverse\rangle )$. Define $\Phi(x,\langle
y,z\rangle )=zg\inverse$. Clearly $\Phi$ is well defined. One can
easily verify that   $\Phi$ over $p$ constitutes a comorphism of
groupoids, which is the global version of the Lie algebroid
comorphism in Example \ref{Ex:morphismGaugealgebroid}.
\end{ex}

\end{document}